\documentclass[11pt,oneside]
{article}

\usepackage{amsfonts}
\usepackage[mathscr]{euscript}
\usepackage{amsthm}
\usepackage{amsmath}
\newcommand{\lbl}[1]{\label{#1}}

\newtheorem{theorem}{Theorem}[section]
\newtheorem{lemma}{Lemma}[section]
\newtheorem{remark}{Remark}[section]
\newtheorem{defi}{Definition}[section]

\newcommand\dd{\displaystyle}
\newcommand\ep{\varepsilon}
\newcommand\ty{\infty}

\newcommand\ud{\underline}

\newcommand{\be}{\begin{equation}}
\newcommand{\ee}{\end{equation}}
\newcommand\bes{\begin{eqnarray}} \newcommand\ees{\end{eqnarray}}
\newcommand{\bess}{\begin{eqnarray*}}
\newcommand{\eess}{\end{eqnarray*}}
\numberwithin{equation}{section}

\begin{document}
\setlength{\baselineskip}{16pt} \pagestyle{myheadings}

\begin{center}{\Large\bf A semilinear parabolic system with a free boundary}\footnote{This work was supported by NSFC Grant 11371113}\\[4mm]
 {\Large  Mingxin Wang,\ \ {Yonggang Zhao}\\[1mm]
{\small Natural Science Research Center, Harbin Institute of Technology, Harbin 150080, PR China}\\[1mm]}
\end{center}

\begin{quote}
\noindent{\bf Abstract}\, This paper deals with a semilinear parabolic system with free boundary in one space dimension. We suppose that unknown functions $u$ and $v$ undergo nonlinear reactions $u^q$ and $v^p$, and exist initially in a interval $\{0\leq x\leq s(0)\}$, but expand to the right with spreading front $\{x=s(t)\}$, with $s(t)$ evolving according to the free boundary condition $s'(t)=-\mu (u_x+\rho v_x)$, where $p,\, q,\, \mu, \,\rho$ are given positive constants. The main purpose of this paper is to understand the existence, uniqueness, regularity and long time behavior of positive solution or maximal positive solution. Firstly, we prove that this problem has a unique positive solution $(u,v,s)$ defined in the maximal  existence interval $[0,T_{\max})$ when $p,\,q\geq 1$, while it has a unique maximal positive solution  $(u,v,s)$ defined in the maximal existence interval $[0,T_{\max})$ when $p<1$ or $q<1$. Moreover,  $(u,v,s)$ and $T_{\max}$ have property that either (i) $T_{\max}=+\ty$, or (ii) $T_{\max}<+\ty$ and
 \bess
  \limsup_{T\nearrow T_{\max}}\|u,\,v\|_{L^{\infty}([0,T]\times[0,s(t)])}=+\infty.
 \eess
Then we study the regularity of $(u,v)$ and $s$. At last, we discuss the global existence ($T_{\max}=+\ty$), finite time blow-up ($T_{\max}<+\ty$), and long time behavior of bounded global solution.


\end{quote}

\section{Introduction and Main Results}
\setcounter{equation}{0}{\setlength\arraycolsep{2pt}

It is well known that free boundary problems for nonlinear parabolic equations have been applied to depict different types of mathematical problems. For instance, it was used in the modeling of ecological dynamics to describe spreading of species \cite{DG, DL,GW, KY, PZ, WMX, WZjdde, WZ},  melting of ice in contact with water \cite{Ru},  chemical vapor deposition in hot wall reactor \cite{MR},  combustion under gravity conditions \cite{GL},  tumor growth \cite{Cui, Tao},  wound healing \cite{FH},  modeling of electrostatic MEMS \cite{ELW, Guoyan}. For rich literatures on free boundary problems and some important theoretical advances, we refer the readers to \cite{CS,Cr,Ru} and the references cited therein.

In this paper, we consider the following semilinear parabolic system with a free boundary
 \bes\label{1.1}
\left\{\begin{array}{ll}
u_t-d_1u_{xx}=v^p, &t>0,\ 0<x<s(t),\\[1mm]
v_t-d_2v_{xx}=u^q, &t>0,\ 0<x<s(t),\\[1mm]
s'(t)=-\mu (u_x+\rho v_x),&t>0,\ x=s(t),\\[1mm]
u_x(t,0)=v_x(t,0)=0,&t>0,\\[1mm]
u(t,s(t))=v(t,s(t))=0,\ \ \ \ &t>0,\\[1mm]
u(0,x)=u_0(x),\ v(0,x)=v_0(x),\ \ &0\leq x\leq s_0,\\[1mm]
s(0)=s_0.
\end{array}\right.
 \ees
Differential equations in (\ref{1.1}) provide a simple example of a reaction diffusion  system. They can be used as a model to describe heat propagation in a two-component  combustible  mixture. In this case $u$ and $v$ represent the temperatures of the interacting components, thermal conductivity is supposed constant for both substances, and heat release  is  described by the power laws.

In problem (\ref{1.1}), $x=s(t)$ represents free boundary which is to be determined together with the solution $(u(t,x),v(t,x))$, parameters $p,\, q,\,d_1,\,d_2,\, \mu, \,\rho$ and $s_0$ are given positive constants, and the assigned initial functions $u_0(x)$ and $v_0(x)$ satisfy
 \bes\label{1.2}
\left\{\begin{array}{ll}
u_0(x),v_0(x)\in W_k^2((0,s_0))\ \ {\rm for \ some} \ k>3,\\[1mm]
 u_0(x), \ v_0(x)>0 \ \ \ {\rm in}\ \ [0,s_0),\\[1mm]
\ u_0'(0)=v_0'(0)=u_0(s_0)=v_0(s_0)=0.
\end{array}\right.
 \ees
Since $k>3$, we have that $W_k^2((0,s_0))\hookrightarrow C^{1+\alpha}([0,s_0])$ with $\alpha=1-1/k$.

Background of the free boundary condition in (\ref{1.1}) can refer to \cite{BDK}. Such kind of free boundary conditions has been used by many authors, please refer to \cite{GW, WMX, WZ} and the references therein.

Many previous mathematical works have been devoted to investigate the corresponding problem on a fixed domain. In particular, Escobedo and Herrero (\cite{EH3}) showed that the problem
 \bes\label{1.3}
\left\{\begin{array}{ll}
u_t- \Delta u=v^p,\ \ &t>0,\ \ x\in\Omega,\\[1mm]
v_t- \Delta v=u^q,&t>0,\ \ x\in\Omega,\\[1mm]
u(t,x)=v(t,x)=0,\ \ \ &t\geq0,\ \ x\in\partial\Omega,\\[1mm]
u(0,x)=u_0(x),v(0,x)=v_0(x),\ \ &x\in\Omega,
\end{array}\right.
 \ees
where $\Omega$ is a bounded domain in $\mathbb{R}^N$ $(N\geq1)$ with smooth boundary, always has a nonnegative unique classical solution provided that either $pq\geq1$, or  $0<pq<1$ and one of the initial functions is different from zero. Moreover, every solution exists for all times if $0<pq\leq1$, but if $pq>1$, solutions may be global or blow up in finite time, according to the size of initial values. When $p>1,\ q>1$, Friedman and Giga (\cite{FG}) established a single point blow up for solutions to (\ref{1.3}) in one space dimension.  In addition, some estimates from above near the blow-up point for a class of positive solutions to (\ref{1.3}) were derived by Caristi and Mitidieri (\cite{CM}) when $\Omega$ is an open ball of $\mathbb{R}^N$ centered at the origin.

When $\Omega=\mathbb{R}^N$, the corresponding problem is the following Cauchy problem
  \bes\label{1.4}
\left\{\begin{array}{ll}
u_t- \Delta u=v^p,\ \ &t>0,\ \ x\in\mathbb{R}^N,\\[1mm]
v_t- \Delta v=u^q,\ \ &t>0,\ \ x\in\mathbb{R}^N,\\[1mm]
u(0,x)=u_0(x),v(0,x)=v_0(x),\ \ &x\in\mathbb{R}^N.
\end{array}\right.
  \ees
A number of properties of solutions to (\ref{1.4}) were acquired in \cite{EH2}. Especially, the solution of (\ref{1.4}) exists globally provided that $0<pq\leq1$. However, if $pq>1$ and
 \bes\label{1.5}
\frac{\kappa+1}{pq-1}\geq \frac N2
 \ees
with $\kappa=\max\{p,q\}$, every nontrivial solution blows up in finite time. On the other hand, if $pq>1$ and (\ref{1.5}) fails, the solution to (\ref{1.4}) might be bounded in any strip $S_T=[0,T)\times \mathbb{R}^N$ or has a finite blow-up time, according to the size of the initial function ($u_0,\,v_0)$. For the case $0<pq<1$, uniqueness result for problem (\ref{1.4}) was established in \cite{EH1}.

If $p=q$ and $u_0=v_0$, then problem (\ref{1.1}) reduces to the following problem
 \bes\label{1.6}
\left\{\begin{array}{ll}
u_t-du_{xx}=u^p, &t>0,\ 0<x<s(t),\\[1mm]
s'(t)=-\mu u_x(t,s(t)), &t>0,\\[1mm]
u_x(t,0)=u(t,s(t))=0,&t>0,\\[1mm]
u(0,x)=u_0(x),&0\leq x\leq s(0)=s_0.
\end{array}\right.
 \ees
When $p>1$, problem (\ref{1.6}) has been studied by Ghidouche et al. (\cite{GST}), Fila and Souplet (\cite{FS}), and Souplet (\cite{S}). The authors of \cite{GST} exhibited an energy condition under which the solution blows up in finite time in $L^\infty$ norm. Moreover, it was shown that all global solutions are bounded and decay uniformly to zero, and that  there are only two possible behaviors for global solutions, either: (i) the solution decays at an exponential rate and the free boundary converges to a finite limit, or (ii) the decay rate of solution is at most polynomial and the free boundary grows up to infinity. In \cite{FS}, it was proved that there exist global solutions with unbounded free boundary and slow decay, i.e. of type (ii). Besides, Souplet (\cite{S}) proved the stability of fast decaying global solution and established a result of continuous dependence of local solution up to the maximum existence time.

If the left fixed boundary $x=0$ in (\ref{1.6}) is replaced by a free boundary $x=r(t)$ governed by $r'(t)=-\mu u_x(t,r(t))$, Zhang and Lin (\cite{ZL}) demonstrated that all  results for (\ref{1.6}) can be extended to the corresponding double free boundary problem.

For simplicity, we introduce the following notations. Assume that $\tau$ is a positive constant and $h(t)\geq\delta>0$ is a continuous function defined in $[0,\tau]$, (or $[0,\tau)$, $(0,\tau]$, $(0,\tau)$). For any given $0\leq\varepsilon<\tau$ and $0\leq\theta\ll 1$, we shall use the following notations, sometimes,
 \bess
 [\varepsilon,\tau]\times[\theta,h(t)]&=&\{(t,x)\in \mathbb{R}^2:\ t\in [\varepsilon,\tau],\ x\in[\theta,h(t)]\},\\[1mm]
 (\varepsilon,\tau]\times[\theta,h(t))&=&\{(t,x)\in \mathbb{R}^2:\ t\in (\varepsilon,\tau],\ x\in[\theta,h(t))\},\\[1mm]
 [\varepsilon,\tau)\times[\theta,h(t))&=&\{(t,x)\in \mathbb{R}^2:\ t\in [\varepsilon,\tau),\ x\in[\theta,h(t))\},\eess
and so on.

For any given $a=(a_1,\cdots,a_n)$ and $b=(b_1,\cdots,b_n)$, we appoint that $a\leq b$ means $a_i\leq b_i$ for all $i$, and $a<b$ means $a_i<b_i$ for all $i$.

\begin{defi}\lbl{d1.1}\, By a positive solution $(u,v,s)$ of $(\ref{1.1})$ defined in $[0,T)$, it means that $s(t)>0$ in $[0,T)$, $u,\,v>0$ in $[0,T)\times[0,s(t))$ and satisfy $(\ref{1.1})$ in the classical sense.

We say that $(u,v,s)$ is a maximal positive solution of $(\ref{1.1})$ defined in $[0,T)$, if for any positive solution $(\hat u,\hat v,\hat s)$ of $(\ref{1.1})$ defined in $[0,\hat T)$ with $\hat T\leq T$, it must hold:
 \[\big(\hat s(t),\,\hat u(t,x),\,\hat v(t,x)\big)\leq\big(s(t),\, u(t,x),\,v(t,x)\big), \ \ \ \forall \ t\in[0,\hat T), \ x\in[0,\hat s(t)].\]
\end{defi}

It should be emphasized that the non-negative and nontrivial solution of (\ref{1.1}) must be positive one since the initial values $u_0(x)>0,\,v_0(x)>0$ in $[0,s_0)$.

Now we state our main results of this paper.

\begin{theorem}\lbl{t1.1}\,There exist a maximum existence time $T_{\max}$ and
 \begin{quote}
{\rm(i)}\, a unique positive solution $(u,v,s)$ of $(\ref{1.1})$ defined in $[0,T_{\max})$ for the case $p\geq 1$ and $q\geq 1$, \\
{\rm (ii)}\, a unique maximal positive solution $(u,v,s)$ of $(\ref{1.1})$ defined in $[0,T_{\max})$ for the case that $p<1$ or $q<1$,
  \end{quote}
such that either $T_{\max}=+\infty$, or $T_{\max}<+\infty$ and
 \bes
  \limsup_{T\nearrow T_{\max}}\|u\|_{L^{\infty}([0,T]\times[0,s(t)])}=+\infty,\ \ \  \limsup_{T\nearrow T_{\max}}\|v\|_{L^{\infty}([0,T]\times[0,s(t)])}=+\infty.
 \lbl{1.7}\ees
\end{theorem}

\begin{remark} When  $p<1$ or $q<1$, it is unfortunately that we can not prove the uniqueness conclusion as {\rm\cite{EH3}} in where the problem {\rm(\ref{1.3})} is concerned.
\end{remark}

\begin{theorem}\lbl{t1.2}\, Let $T>0$ and $(u,v,s)$ be a positive solution of $(\ref{1.1})$ defined in $[0,T]$. Then
 \bess
 u,\,v\in C^{1+\frac\beta 2,\,2+\beta}(D_T),
  \ \ \ s\in C^{1+\frac{1+\beta}2}([0,T]),
  \eess
where $D_T=(0,T]\times[0,s(t)]$, $\beta=\min\{p,\,q\}$. Particularly,

{\rm(i)}\, if $p$ and $q$ are positive integers, then $u,\,v\in C^\infty(D_T)$, $s\in C^{\ty}((0,T])$;

{\rm(ii)}\, if $p,q\geq 1$ and one of them is not integer, then the regularity of $(u,v,s)$ depends strongly on the relationship between parameters $p,\,q$ and $\alpha:=1-3/k$. Here we only give the result for a special case. Take $q\geq p=1+\lambda$ and $0<\lambda\leq\alpha/2$, then
 \bess u,\,v\in C^{1+\frac{1+\lambda}2,\,3+\lambda}(D_T), \ \
 s\in C^{2+\frac{\lambda}2}((0,T]).
  \eess
\end{theorem}

\begin{theorem}\lbl{t1.3}
Let $(u,v,s)$ be a positive solution of problem $(\ref{1.1})$ defined in $[0,T]$, then $s'(t)>0$ in $(0,T]$.
\end{theorem}

\begin{theorem}\label{t1.4} \, Let $s_0$, $\mu$ and $\rho$ be fixed, $(u,v,s)$ and $T_{\max}$ be obtained in Theorem $\ref{t1.1}$.

{\rm(i)} Assume that $pq>1$. Then $T_{\max}=+\ty$, i.e., $(u,v)$ exists globally in time provided that initial functions $u_0(x)$ and $v_0(x)$ are suitably small; while $T_{\max}<+\ty$, i.e., $(u,v)$ will blow up in finite time provided that initial functions $u_0(x)$ and $v_0(x)$ are large enough.

{\rm(ii)} If $pq\leq1$, then $T_{\max}=+\ty$.
\end{theorem}

\begin{theorem}\lbl{t1.5}\, Let $s_0$, $\mu$ and $\rho$ be fixed, $(u,v,s)$ and $T_{\max}$ be obtained in Theorem $\ref{t1.1}$. If $T_{\max}=+\ty$, $s_{\infty}:=\lim_{t\to+\infty}s(t)<+\infty$, $u$ and $v$ are bounded, then
 \bess
 \lim_{t\to+\infty}\,\max_{0\leq x\leq s(t)}u(t,x)=\lim_{t\to+\infty}\,\max_{0\leq x\leq s(t)}v(t,x)=0.
 \eess
\end{theorem}

The plan of this article is as follows. We first prove Theorems \ref{t1.2} and \ref{t1.3} in section \ref{s.2} rather than Theorem \ref{t1.1} because the proof of Theorem \ref{t1.1} is very complicated. Section \ref{s.3} is devoted to deal with Theorem \ref{t1.1} for the case $p,\,q\geq 1$. In section \ref{s.4}, we establish two comparison principles which will be used in the last two sections. The proof of Theorem \ref{t1.1} for the case that either $p<1$ or $q<1$ will be given in section \ref{s.5}. In the last section, we shall deal with Theorems \ref{t1.4} and \ref{t1.5}.

Since parameters $d_1,\,d_2,\,p,\,q,\,\mu$ and $\rho$ are fixed, we don't emphasize the dependence of the generic estimated constants on these parameters at each step of the following estimates.

\section{Proofs of Theorems \ref{t1.2} and \ref{t1.3}}\label{s.2}

In this section we first study the regularity of positive solution $(u,v,s)$, and then present the monotonicity of the free boundary $s(t)$. That is, Theorems \ref{t1.2} and \ref{t1.3} are proved successively.

In order to show Theorem \ref{t1.2}, we transform the free boundary problem (\ref{1.1}) into an initial-boundary value problem with fixed boundary. And then applying the Schauder interior estimate for parabolic equations, we get the regularity of $(u,v,s)$. However, since the regularity of reaction terms heavily depends on the values of $p$ and $q$, we have to divide them into three cases: (i) $p$ and $q$ are positive integers;  (ii) $p,q\geq 1$ and one of them is not integer; (iii) $p<1$ or $q<1$.

{\bf Proof of Theorem \ref{t1.2}} \, Note that $u_0,v_0$ satisfy (\ref{1.2}). Applying the $L^p$ theory for parabolic equations and Sobolev embedding theorem, and then combining with the free boundary condition
 $s'(t)=-\mu\big[u_x(t,s(t))+\rho v_x(t,s(t))\big],$
it is not difficult to derive (see, for example, \cite{GW, WZ}) that
  \bes
  u,\,v\in C^{\frac{1+\alpha} 2,\,1+\alpha}(D_T), \ \ s\in C^{1+\frac\alpha 2}([0,T])\lbl{2.1}\ees
with $\alpha=1-3/k$. Clearly, $u_x(t,s(t))\leq 0$, $v_x(t,s(t))\leq 0$ since $u,v>0$ in $(0,T]\times(0,s(t))$ and $u,v=0$ at free boundary $x=s(t)$. Hence $s'(t)\geq 0$.

The idea of the following proof comes from \cite{Sc}. We shall use the transformation
  \bes
 y=x/s(t),\  \ \ w(t,y)=u(t,x),\ \ \ z(t,y)=v(t,x)
 \lbl{2.2}\ees
to straighten the free boundary $x=s(t)$. A series of detailed calculation asserts
 \bes
\left\{\begin{array}{lll}
u_t=w_t-ys'(t)s^{-1}(t)w_y,\ \ v_t=z_t-ys'(t)s^{-1}(t)z_y,\\[1mm]
u_x=s^{-1}(t)w_y, \ \ \ \ v_x=s^{-1}(t)z_y, \\[1mm]
u_{xx}=w_{yy}s^{-2}(t),\ \ v_{xx}=z_{yy}s^{-2}(t),
\end{array}\right.\lbl{2.3}
 \ees
and
 \bes\label{2.4}
\left\{\begin{array}{ll}
w_t-d_1f(t)w_{yy}-g(t,y)w_y=z^p, &0<t\leq T,\ 0<y<1,\\[1mm]
z_t-d_2f(t)z_{yy}-g(t,y)z_y=w^q, &0<t\leq T,\ 0<y<1,\\[1mm]
w_y(t,0)=z_y(t,0)=w(t,1)=z(t,1)=0, \ \ &0<t\leq T,\\[1mm]
w(0,y)=u_0(s_0y),\ z(0,y)=v_0(s_0y),&0\leq y\leq 1,
\end{array}\right.
 \ees
where $f(t)=s^{-2}(t)$, $g(t,y)=ys'(t)/s(t)$. This is an initial-boundary value problem with fixed boundary.

We first consider the case that $p,\,q\geq 1$. By (\ref{2.1}), it can be deduced that $v^p,\,u^q\in C^{\frac\alpha2,\alpha}(D_T)$ and $s'\in C^{\frac{\alpha} 2}((0,T])$.
And so
 \[z^p,\,w^q\in C^{\frac\alpha2,\alpha}\big((0,T]\times[0,1]\big), \ \
 f\in C^{1+\frac\alpha 2}((0,T]), \ \ g\in C^{\frac{\alpha}2,\alpha}\big((0,T]\times[0,1]\big)\]
since $s(t)\geq s_0$. For any given $0<\ep\ll 1$, applying Theorem 10.1 of \cite[Chap.4, p.351]{LSU} to problem (\ref{2.4}) in $[\varepsilon,T]\times[\ep,1]$ and $[\varepsilon,T]\times[0,1-\ep]$, respectively, we obtain that
  \[w,\,z\in C^{1+\frac{\alpha}2,\,2+\alpha}\big([\varepsilon,T]\times[\ep,1]\big)
  \bigcap C^{1+\frac{\alpha}2,\,2+\alpha}\big([\varepsilon,T]\times[0,1-\ep]\big).\]
Taking advantage of (\ref{2.2}) and (\ref{2.3}), it follows that
 \bess u,\,v\in C^{1+\frac{\alpha}2,\,2+\alpha}\big([\varepsilon,T]\times[\ep s(t),s(t)]\big)
  \bigcap C^{1+\frac{\alpha}2,\,2+\alpha}\big([\varepsilon,T]\times[0,(1-\ep)s(t)]\big).
  \eess
Due to the arbitrariness of $\ep$, one achieves
 \[u,\,v\in C^{1+\frac{\alpha}2,\,2+\alpha}(D_T),  \ \ \Longrightarrow \ u_x,\,v_x\in C^{\frac{1+\alpha}2,\,1+\alpha}(D_T).\]
Hence, by the condition $s'(t)=-\mu\big[u_x(t,s(t))+\rho v_x(t,s(t))\big]$, it is immediately to get $s'\in C^{\frac{1+\alpha}2}((0,T])$.

\vskip 4pt (i)\, When both $p$ and $q$ are positive integers, we still have $u^q,\,v^p\in C^{1+\frac{\alpha} 2,\,2+\alpha}(D_T)$ since $u,\,v\in C^{1+\frac{\alpha}2,\,2+\alpha}(D_T)$. Consequently,
\[z^p,\,w^q\in C^{1+\frac{\alpha}2,2+\alpha}\big((0,T]\times[0,1]\big), \ \
 f\in C^{1+\frac{1+\alpha} 2}((0,T]), \ \ g\in C^{\frac{1+\alpha}2,1+\alpha}\big((0,T]\times[0,1]\big).\]
Let $\alpha_1=1+\alpha$. Similar to the above, it can be deduced that
 \bess u,\,v\in C^{1+\frac{\alpha_1} 2,\,2+\alpha_1}(D_T)=C^{1+\frac{1+\alpha} 2,\,3+\alpha}(D_T), \ \
 s'\in C^{\frac{1+\alpha_1}2}((0,T])=C^{1+\frac{\alpha}2}((0,T]).
  \eess
Repeating such processes, the desired result will be obtained eventually.

(ii)\, When $q\geq p=1+\lambda$ and $0<\lambda\leq\alpha/2$, the derivative $(v^p)_{xx}$ does not exist at the point $(t,s(t))$ no matter how smooth $v$ is, since $v(t,s(t))=0$. Even though we have known that $v\in C^{1+\frac{\alpha}2,\,2+\alpha}(D_T)$,  it can be only obtained that  $v^\lambda\in C^\lambda$ in $t$ and $x$ due to $0<\lambda<1$, and $v_t\in C^{\frac\alpha 2,\alpha}(D_T)$, $v_x\in C^{\frac{1+\alpha}2,\,1+\alpha}(D_T)$. Notice $\lambda\leq\alpha/2$, one can only get $(v^p)_t=pv^\lambda v_t\in C^\lambda$, $(v^p)_x=pv^\lambda v_x\in C^\lambda$ in $t$ and $x$. Hence, $v^p\in C^{1+\lambda}$ in $t$ and $x$. Using the standard notation,
$v^p\in C^{\frac{1+\lambda}2,\,1+\lambda}(D_T)$. Certainly, $u^q\in C^{\frac{1+\lambda}2,\,1+\lambda}(D_T)$ since $q\geq p=1+\lambda$. Note $s'\in C^{\frac{1+\alpha} 2}((0,T])$ and $1+\lambda<1+\alpha$. In the same way as the argument of case (i), it can be deduced that
  \bes u,\,v\in C^{1+\frac{1+\lambda}2,\,3+\lambda}\big((0,T]\times[0,s(t)]\big), \ \
 s'\in C^{1+\frac{\lambda}2}((0,T]).\lbl{2.5}
  \ees
Since $\lambda<(1+\lambda)/2$, similar to the above, we can only deduce $v^p\in C^{\frac{1+\lambda}2,\,1+\lambda}(D_T)$ and (\ref{2.5}), but cannot arrive at the higher regularity of $(u,v,s)$.

At last, we consider the situation that at least one of the exponents $p,\ q$ is less than one. For convenience, we only sketch how to deal with the case where $q\geq 1>p$, and leave the general situation to reader. Since $v(t,s(t))=0$, the derivative $(v^p)_x$ does not exist at $(t,s(t))$  no matter how smooth  $v$ is. So one can only obtain that, by (\ref{2.1}), $v^p\in C^{p\frac{1+\alpha}2}$ in $t$ and $v^p\in C^{p}$ in $x$. Using the standard notation, $v^p\in C^{p/2,\,p}(D_T)$. Same as the argument of (ii), we acquire
 \[u,\,v\in C^{1+\frac p2,\,2+p}(D_T), \ \ s'\in C^{\frac{1+p}2}((0,T]),\]
but cannot get the higher regularity of $(u,v,s)$. The proof is finished. \ \ \ \fbox{}

\vskip 4pt {\bf Proof of Theorem \ref{t1.3}} \, When $p,\,q\geq 1$, we know that $s\in C^2((0,T])$ by Theorem \ref{t1.2}. The domain $D_T$ has an {\it interior sphere} property at the right boundary $x=s(t)$.
Using Hopf's boundary lemma for parabolic equations to the first and second equations of (\ref{1.1}) yields
  \bes
  u_x(t,s(t))<0,\ \ v_x(t,s(t))<0,\ \ {\rm for}\ \ 0<t\leq T.
  \lbl{2.6}\ees
And then substituting the above two inequalities into the third equation of (\ref{1.1}) leads to the desired result $s'(t)>0$ in $(0,T]$.

If $p<1$ or $q<1$, we only know $s\in C^{1+\frac{1+\beta} 2}([0,T])$ by Theorem \ref{t1.2}, and then cannot guarantee that the domain $D_T$ has an {\it interior sphere} property at the right boundary $x=s(t)$. Hence, the Hopf boundary lemma cannot be used directly to (\ref{1.1}). To overcome this, we use the transformation (\ref{2.2}) to straighten the free boundary $x=s(t)$. Since $w,z>0$ in $[0,T]\times[0,1)$ and $w(t,1)=z(t,1)=0$, applying the Hopf boundary lemma to (\ref{2.4}) we get $w_y(t,1)<0$, $z_y(t,1)<0$. By virtue of $u_x=s^{-1}(t)w_y$, $v_x=s^{-1}(t)z_y$, (\ref{2.6}) is deduced. Therefore, $s'(t)>0$ in $(0,T]$. The proof is complete. \ \ \ \ \fbox{}

\section{Proof of Theorem \ref{t1.1} for the case $p,\,q\geq 1$}\label{s.3}

In this section, we shall give the proof of Theorem \ref{t1.1} for the case $p,\,q\geq 1$ by means of the contraction mapping theorem and extension method. That is, we shall prove the following theorem

\begin{theorem}\lbl{t3.1} \, Let $p\geq 1$ and $q\geq 1$. Then there exist a maximum existence time $T_{\max}$ and
a unique positive solution $(u,v,s)$ of $(\ref{1.1})$ defined in $[0,T_{\max})$,
such that either $T_{\max}=+\infty$, or $T_{\max}<+\infty$ and $(\ref{1.7})$ holds.
 \end{theorem}

The proof of Theorem \ref{t3.1} is based on the following three lemmas. In the first one we show that (\ref{1.1}) has a unique local solution which can be extended to the maximal existence interval $(0,\,T_{\max})$. In the second lemma we give the estimate of upper bound of $s'(t)$. In the last lemma we prove that (\ref{1.7}) holds when $T_{\rm max}<\infty$.

\begin{lemma}\lbl{l3.1} \, If $p\geq 1$ and $q\geq 1$, Then there exist a maximum existence time $T_{\max}$ and a unique positive solution $(u,v,s)$ of $(\ref{1.1})$ defined in $[0,T_{\max})$.
\end{lemma}

{\bf Proof.}\, This proof is divided into two steps. In the first step, we introduce the standard transformation to straighten the free boundary, then take advantage of the contraction mapping theorem to show the local existence and uniqueness. In the second step, we extend the unique solution to the maximal existence interval $(0,\,T_{\max})$.

{\it Step 1}\, Denote $\hat s=-\mu(u_0'(s_0)+\rho v_0'(s_0))$. We first prove that there exists $0<T\ll 1$, depending only on $s_0$, $\hat s$, $\|u_0\|_{W_k^2((0,s_0))}$ and $\|v_0\|_{W_k^2((0,s_0))}$, such that problem $(\ref{1.1})$ admits a unique positive solution $(u,v,s)$ defined in $[0,T]$.

\vskip 4pt This proof can be done by modifying the arguments of \cite{CF, DL, GW, WZ}. We provide the details here for the reader¡¯s convenience. Let $\zeta(y)$ be a function in $C^3[0,+\infty)$ satisfying
 \begin{center}
$\zeta(y)=1$ \ if \ $\dd|y-s_0|<\frac{s_0}4$,\ \ $\zeta(y)=0$ \ if \ $|y-s_0|>\dd\frac{s_0}2$,\ \ $|\zeta'(y)|<\dd\frac{6}{s_0}$ \ for all \ $y$.
 \end{center}
Define
$$(t,x)\rightarrow(t,y),\ {\rm where}\ \  x=y+\zeta(y)(s(t)-s_0),\ 0\leq y<+\infty.$$
Note that for fixed $t>0$, as long as
   $$|s(t)-s_0|\leq\frac{s_0}8,$$
the transformation $(t,x)\rightarrow(t,y)$ is a diffeomorphism  from $[0,+\infty)$ onto $[0,+\infty)$. Moreover,
 \bess
\left.\begin{array}{rcl}
x=s(t)&\Longleftrightarrow& y=s_0,\\[1mm]
0\leq x\leq s(t)&\Longleftrightarrow& 0\leq y\leq s_0.
\end{array}\right.
 \eess
Direct calculations indicate
\bes
\left\{\begin{array}{ll}
\displaystyle{\frac{\partial y}{\partial x}}=\frac 1{1+\zeta'(y)(s(t)-s_0)}\equiv \sqrt{A(y,s(t))},\\[4mm]
\displaystyle\frac{\partial^2 y}{\partial x^2}=\frac{-\zeta''(y)(s(t)-s_0)}{ [1+\zeta'(y)(s(t)-s_0)]^3}\equiv {B(y,s(t))},\\[4mm]
\displaystyle\frac{\partial y}{\partial t}=\frac{-s'(t)\zeta(y)}{1+\zeta'(y)(s(t)-s_0)}\equiv -s'(t)C(y,s(t)).
\end{array}\right.\lbl{3.2}
\ees
If we set
\bes
\left\{\begin{array}{c}
u(t,x)=u(t,y+\zeta(y)(s(t)-s_0))=U(t,y),\\[1.5mm]
v(t,x)=v(t,y+\zeta(y)(s(t)-s_0))=V(t,y),
\end{array}\right.\lbl{3.3}
\ees
then $(U,V,s)$ satisfies
 \bes\label{3.4}
\left\{\begin{array}{ll}
U_t-d_1AU_{yy}-(d_1B+s'(t)C)U_y=V^p, \ \ &t>0,\ 0<y<s_0,\\[1mm]
V_t-d_2AV_{yy}-(d_2B+s'(t)C)V_y=U^q, &t>0,\ 0<y<s_0,\\[1mm]
s'(t)=-\mu (U_y+\rho V_y),&t>0,\ y=s_0,\\[1mm]
U_y(t,0)=V_y(t,0)=0,&t>0,\\[1mm]
U(t,s_0)=V(t,s_0)=0,\ \ \ \ &t>0,\\[1mm]
U(0,y)=U_0(y),\ V(0,y)=V_0(y),&0\leq y\leq s_0,
\end{array}\right.
 \ees
where $A=A(y,s(t))$, $B=B(y,s(t))$, $C=C(y,s(t))$, $U_0(y)=u_0(y)$, $V_0(y)=v_0(y)$.

\vskip 4pt Obviously, $\hat s=-\mu(u_0'(s_0)+\rho v_0'(s_0))\geq 0$. For $0<T\leq\displaystyle\frac{s_0}{8(1+\hat s)}$, define
 \bess
 Q_T&=&[0,T]\times[0,s_0],\\[1mm]
\mathcal{X}_{1T}&=&\{U\in C(Q_T):\ U\geq0,\ U(0,y)=u_0(y),\ \|U-u_0\|_{C(Q_T)}\leq1\},\\[1.5mm]
\mathcal{X}_{2T}&=&\{V\in C(Q_T):\ V\geq0,\ V(0,y)=v_0(y),\ \|V-v_0\|_{C(Q_T)}\leq1\},\\[1.5mm]
\mathcal{X}_{3T}&=&\{s\in C^1([0,T]):\ s(0)=s_0,\  s'(0)=\hat s,\ \|s'-\hat s\|_{C([0,T])}\leq1\}.
 \eess
It is not difficult to verify that $\mathcal{X}_{T}:=\mathcal{X}_{1T}\times\mathcal{X}_{2T}\times\mathcal{X}_{3T}$ is a closed convex set in $C(Q_T)\times C(Q_T)\times C^1([0,T])$.

Next, we shall prove the existence result by means of the contraction mapping theorem. Firstly, it is easy to show that, for arbitrary $(U,V,s)\in \mathcal {X}_T$,
  $$|s(t)-s_0|\leq T(1+\hat s)\leq \frac{s_0}8$$
since $T\leq\displaystyle\frac{s_0}{8(1+\hat s)}$. Therefore, the above transformation $(t,x)\rightarrow(t,y)$ is well defined.

\vskip 4pt For any $(U,V,s)\in \mathcal {X}_T$, we consider the following initial boundary value problem
 \bes\label{3.5}
\left\{\begin{array}{ll}
\tilde U_t-d_1A\tilde U_{yy}-(d_1B+s'(t)C)\tilde U_y=V^p, &t>0,\ 0<y<s_0,\\[1mm]
\tilde U_y(t,0)=0,\ \tilde U(t,s_0)=0,\ \ \ \ &t>0,\\[1mm]
\tilde U(0,y)=U_0(y),&0\leq y\leq s_0.
\end{array}\right.
 \ees
Since $U_0\in W^2_k((0,s_0))$, the standard partial differential equations theory \cite{Fr,LSU} illustrates that problem (\ref{3.5}) admits a unique solution $\tilde U\in C^{\frac{1+\alpha}2,1+\alpha}(Q_T)$ with
 \bes\label{3.6}
\|\tilde U\|_ {C^{\frac{1+\alpha}2,1+\alpha}(Q_T)}\leq C_1,  \ \ \alpha=1-3/k,
 \ees
where $C_1$ is a positive constant depending on $s_0$, $\hat s$, $\|u_0\|_{W_k^2((0,s_0))}$ and $\|v_0\|_{W_k^2((0,s_0))}$. Similarly, for given $(U,V,s)\in \mathcal {X}_T$, initial-boundary value problem
 \bess
\left\{\begin{array}{ll}
\tilde V_t-d_2A\tilde V_{yy}-(d_2B+s'(t)C)\tilde V_y=U^q, &t>0,\ 0<y<s_0,\\[1mm]
\tilde V_y(t,0)=0,\ \tilde V(t,s_0)=0,\ \ \ \ &t>0,\\[1mm]
\tilde V(0,y)=V_0(y),&0\leq y\leq s_0
\end{array}\right.
 \eess
has a unique solution $\tilde V\in C^{\frac{1+\alpha}2,1+\alpha}(Q_T)$ and
 \bes\label{3.7}
 \|\tilde V\|_ {C^{\frac{1+\alpha}2,1+\alpha}(Q_T)}\leq C_1.
 \ees

Define
$$\tilde s(t)=s_0-\mu\int_0^t(\tilde U_y(\tau,s_0)+\rho {\tilde V}_y(\tau,s_0))d\tau,$$
then  $$\tilde s(0)=s_0,\ \tilde s'(0)=\hat s,\ \tilde s'(t)=-\mu(\tilde U_y+\rho {\tilde V}_y)(t,s_0).$$
Therefore, $\tilde s'(t)\in C^{\frac\alpha 2}([0,T])$ and
 \bes\label{3.8}
  \|\tilde s'(t)\|_{C^{\frac\alpha 2}([0,T])}\leq \mu(1+\rho)C_1:=C_2.
 \ees

Now, we introduce a mapping $\mathcal{F}:\ \mathcal{X}_T\rightarrow C(Q_T)\times C(Q_T)\times C^1([0,T])$ by
$$\mathcal{F}(U,V,s)=(\tilde U,\tilde V,\tilde s).$$
We next prove that $\mathcal{F}$ has a unique fixed point, which is a solution to system (\ref{3.4}). In view of (\ref{3.6})--(\ref{3.8}), we have that
 \bess
\left.\begin{array}{l}
\|\tilde U-u_0\|_{C(Q_T)}\leq\|\tilde U\|_{C^{\frac{1+\alpha} 2,0}(Q_T)}T^{\frac{1+\alpha} 2}\leq C_1T^{\frac{1+\alpha} 2},\\[2.5mm]
\|\tilde V-v_0\|_{C(Q_T)}\leq\|\tilde V\|_{C^{\frac{1+\alpha}2,0}(Q_T)}T^{\frac{1+\alpha} 2}\leq C_1T^{\frac{1+\alpha}2},\\[2.5mm]
\|\tilde s'-\hat s\|_{C([0,T])}\leq\|\tilde s'\|_{C^{\frac\alpha 2}([0,T])}T^{\frac\alpha 2}\leq C_2T^{\frac\alpha 2}.
\end{array}\right.
 \eess
So, if we choose
  $$T\leq T_0:=\min\left\{C_1^{-\frac 2{1+\alpha}},\  C_2^{-\frac 2\alpha}, \ \frac{s_0}{8(1+\hat s)}\right\},$$
then $\mathcal{F}$ maps $\mathcal{X}_T$ into itself.

Now, it will be showed that $\mathcal{F}$ is a contraction mapping on $\mathcal{X}_T$ for sufficiently small $T>0$. In fact, let $(U_i,V_i,s_i)\in \mathcal{X}_T$ for $i=1,2$ and denote $(\tilde U_i,\tilde V_i,\tilde s_i)=\mathcal{F}(U_i,V_i,s_i)$. By virtue of (\ref{3.6})--(\ref{3.8}), it is easy to obtain
 $$\|\tilde U_i\|_{C^{\frac{1+\alpha}2,1+\alpha}(Q_T)}\leq C_1,\ \ \|\tilde V_i\|_{C^{\frac{1+\alpha}2,1+\alpha}(Q_T)}\leq C_1,\ \ \|\tilde s'_i\|_{C^{\frac\alpha 2}([0,T])}\leq C_2.$$
Set $\hat U=\tilde U_1-\tilde U_2$ and $\hat V=\tilde V_1-\tilde V_2$, it can be verified that $\hat U$ satisfies
 \bess
\left\{\begin{array}{ll}
\hat U_t-d_1A(y,s_2)\hat U_{yy}-(d_1B(y,s_2)+s'(t)C(y,s_2))\hat U_y=\mathcal{U}, &t>0,\ 0<y<s_0,\\[1mm]
\hat U_y(t,0)=0,\ \hat U(t,s_0)=0,\ \ \ \ &t>0,\\[1mm]
\hat U(0,y)=U_0(y),&0\leq y\leq s_0,
\end{array}\right.
  \eess
where
 \begin{eqnarray*}
  \mathcal{U}=& d_1[A(y,s_1)-A(y,s_2)]\tilde U_{1,yy}+d_1[B(y,s_1)-B(y,s_2)]\tilde U_{1,y} \\[1mm]
   &+[s_1'C(y,s_1)-s_2'C(y,s_2)]\tilde U_{1,y}+V_1^p-V_2^p.
 \end{eqnarray*}
Again, applying the $L^p$ estimates for parabolic equations and Sobolev's imbedding theorem, it is deduced that
 \bes\label{3.9}
 \|\hat U\|_{C^{\frac{1+\alpha}2,{1+\alpha}}(Q_T)}\leq C_3\left(\|V_1-V_2\|_{C(Q_T)}+\|s_1-s_2\|_{C^1([0,T])}\right),
 \ees
for some positive constant $C_3$ which depends on the $L^\infty$-norms of functions $A,B,C$, and constants $C_1$ and $C_2$. Similarly,
 \bes\label{3.10}
 \|\hat V\|_{C^{\frac{1+\alpha}2,{1+\alpha}}(Q_T)}\leq C_3\left(\|U_1-U_2\|_{C(Q_T)}+\|s_1-s_2\|_{C^1([0,T])}\right).
 \ees
Taking the difference of equations for $\tilde s_1,\ \tilde s_2$ leads to
 \bes\label{3.11}
 \|\tilde s_1'-\tilde s_2'\|_{C^{\frac\alpha 2}([0,T])}\leq C_4\left(\|\hat U_y\|_{C^{\frac\alpha 2,0}(Q_T)}+\|\hat V_y\|_{C^{\frac\alpha 2,0}(Q_T)}\right),
 \ees
where $C_4=\mu(1+\rho)C_3$. Using (\ref{3.9})--(\ref{3.11}), and assuming that $T\leq 1$, we obtain
 \begin{eqnarray*}
  &\|\hat U\|_{C^{\frac{1+\alpha}2,{1+\alpha}}(Q_T)}+ \|\hat V\|_{C^{\frac{1+\alpha}2,{1+\alpha}}(Q_T)}+\|\tilde s_1'-\tilde s_2'\|_{C^{\frac\alpha 2}([0,T])}\\[1mm]
  &\leq C_5\left(\|U_1-U_2\|_{C(Q_T)}+\|V_1-V_2\|_{C(Q_T)}+\|s_1'-s_2'\|_{C([0,T])}\right)
 \end{eqnarray*}
for some positive constant $C_5$ depending only on $C_3$ and $C_4$, where we have used the facts that $s_1(0)=s_2(0)$ and $\|s_1-s_2\| _{C([0,T])} \leq T \|s_1'-s_2'\|_{C([0,T])}$. Hence, if we select
  $$T=\min\left\{T_0, \ (2C_5)^{-\frac 2\alpha}\right\}=\min\left\{C_1^{-\frac 2{1+\alpha}}, \ C_2^{-\frac 2\alpha}, \ \frac{s_0}{8(1+\hat s)}, \ (2C_5)^{-\frac 2\alpha}\right\},$$
then
 \begin{eqnarray*}
  &&\|\hat U\|_{C(Q_T)}+ \|\hat V\|_{C(Q_T)}+\|\tilde s_1'-\tilde s_2'\|_{C([0,T])}\\[1mm]
  &&\leq T^{\frac\alpha 2}\left(\|\hat U\|_{C^{\frac{1+\alpha}2,{1+\alpha}}(Q_T)}+ \|\hat V\|_{C^{\frac{1+\alpha}2,{1+\alpha}}(Q_T)}+\|\tilde s_1'-\tilde s_2'\|_{C^{\frac\alpha 2}([0,T])}\right)\\[1mm]
  &&\leq C_5T^{\frac\alpha 2}\left(\|U_1-U_2\|_{C(Q_T)}+\|V_1-V_2\|_{C(Q_T)}+\|s_1'-s_2'\|_{C([0,T])}\right)\\[1mm]
  &&\leq\frac 12\left(\|U_1-U_2\|_{C(Q_T)}+\|V_1-V_2\|_{C(Q_T)}+\|s_1'-s_2'\|_{C([0,T])}\right),
\end{eqnarray*}
and  $\mathcal{F}$ maps $\mathcal{X}_T$ into itself.

The above arguments demonstrate that $\mathcal{F}$ is a contraction mapping on $\mathcal{X}_T$ if $T$ is small enough. Thus we can apply the contraction mapping theorem to conclude that $\mathcal{F}$ has a unique fixed point $(U,V,s)$ in $\mathcal{X}_T$. From the preceding discussions we also see that $U,\,V\in C^{\frac{1+\alpha}2,1+\alpha}(Q_T)$, $s\in C^{1+\frac\alpha 2}([0,T])$ and the corresponding (\ref{3.6})--(\ref{3.8}) hold. This shows that $(U,V,s)$ is the unique weak solution of (\ref{3.4}). By virtue of $s\in C^{1+\frac\alpha 2}([0,T])$ and property of $\zeta(y)$, it can be seen that coefficients of (\ref{3.4}) belong to $C^{\frac\alpha 2,\,\alpha}(Q_T)$. Applying Schauder's estimate in $[\varepsilon,T]\times[0,s_0]$ for any $0<\varepsilon<T$, it is  derived that
  \bes
  U,\,V\in C^{1+\frac\alpha 2,2+\alpha}\big((0,T]\times[0,s_0]\big),\lbl{3.12}\ees
which implies that $(U,V,s)$ is the unique classical solution of problem (\ref{3.4}) defined in $[0,T]$.

Recalling (\ref{3.2}), (\ref{3.3}) and $s\in C^{1+\frac\alpha 2}([0,T])$, in view of the properties of function $\zeta(y)$ and transformation
$(t,x)\rightarrow(t,y)$, it follows from (\ref{3.12}) that
 \[u,\,v\in C^{1+\frac\alpha 2,2+\alpha}((0,T]\times[0,s(t)]), \ \ \alpha=1-3/k.\]
So, by Sobolev's imbedding theorem, $s'\in C^{\frac{1+\alpha}2}((0,T])$. This suggests that $(u,v,s)$ is the unique positive solution of (\ref{1.1}) defined in $[0,T]$.

\vskip 4pt {\it Step 2} \, Since the uniqueness result holds, the solution $(u,v,s)$ can be extended to $[T,T+\delta]$ for some $\delta>0$ using the above method. Repeating this procedure, we can define
 \[T_{\max}=\sup\big\{T>0:(u,v,s) \ {\rm is \ the \ unique \ positive \ solution \ of \ (\ref{1.1}) \ defined \ in }\ [0,T]\big\}.\]
Then $(u,v,s)$ is the unique positive solution of (\ref{1.1}) defined in $[0,T_{\max})$ and satisfies
  \[u,\,v\in C^{1+\frac\alpha 2,2+\alpha}\big((0,T_{\max})\times[0,s(t)]\big),
  \ \ \ s'\in C^{\frac{1+\alpha}2}(0,T_{\max}).\]
The proof is complete. \ \ \ \fbox{}

In order to prove (\ref{1.7}), we first give an estimate of $s'(t)$ when $u$ and $v$ are bounded.

\begin{lemma}\label{l3.2}\, Suppose that $p,q>0$ and $(u,v,s)$ is a positive solution of $(\ref{1.1})$ defined in $[0,T)$ for some $T\in (0,+\infty)$. If $u$ and $v$ are bounded for $t\in [0,T)$ and $x\in[0,s(t)]$, then there exists a positive constant $C$ independent of $T$ such that $0<s'(t)\leq C$ for $t\in(0,T)$.
\end{lemma}

{\bf Proof.}\,  We have known that $s'(t)>0$ in $(0,T)$ by Theorem \ref{t1.2}. Since $u$ and $v$ are bounded for $t\in [0,T)$ and $x\in[0,s(t)]$, it is easy to see that $s(t)$ is bounded in $[0,T)$. Let $M$ be the bound of $u$ and $v$. We shall compare $u$ and $v$ with some auxiliary functions (see \cite{DL} or \cite{DG}). To do this, define a comparison function by
  $$w(t,x)=M\left[2K(s(t)-x)-K^2(s(t)-x)^2\right]$$
for some appropriate positive constant $K$ over region
$$D^K=\{(t,x):\ 0<t<T,\ s(t)-1/K<x<s(t)\}.$$

First of all, one can easily compute that, for any $(t,x)\in D^K$,
  \begin{eqnarray*}
  &&w_t=2MK[1-K(s(t)-x)]s'(t)\geq0,\\[1mm]
  &&-w_{xx}=2MK^2,\ \ \ v^p\leq M^p.
\end{eqnarray*}
It follows that, if $K^2\geq \frac{M^{p-1}}{2 d_1}$, then
 \[w_t-d_1w_{xx}\geq2d_1MK^2\geq v^p=u_t-d_1u_{xx}\ \ \ {\rm in}\ \ D^K.\]
On the other hand, it is clear that
  $$w(t,s(t)-{K}^{-1})=M\geq u(t,s(t)-{K}^{-1}),\ \ w(t,s(t))=u(t,s(t))\ \ \ {\rm in}\ \ (0,T).$$
As long as $K$ is further chosen such that
 \bes\label{3.13}
 u_0(x)\leq w(0,x)\ \ \ {\rm in}\ \ [s_0-K^{-1},s_0],
 \ees
then $u(t,x)\leq w(t,x)$ for $(t,x)\in D^K$ by use of the maximum principle to $w-u$ over $D^K$. And it then follows  that
 \bes\label{3.14}
  u_x(t,s(t))\geq w_x(t,s(t))=-2MK.
 \ees

Now we prove that there exists $K$ independent of $T$ such that (\ref{3.13}) holds. Direct calculation gives
$$w_x(0,x)=-2MK[1-K(s_0-x)]\leq -MK\ \ \ {\rm on}\ \ [s_0-(2K)^{-1},\, s_0].$$
Hence, for
 $$K=\max\left\{\frac{4\|u_0\|_{C^1([0,s_0])}}{3M},\ \ \left(\frac{M^{p-1}}{2d_1}\right)^{1/2}\right\},$$
there holds
 $$w_x(0,x)\leq-{\frac{4\|u_0\|_{C^1([0,s_0])}}3}\leq u'_0(x)\ \ \ {\rm on}\ \ [s_0-(2K)^{-1},\, s_0],$$
and then integrating the above inequality over $[x,s_0]$ and using $w(0,s_0)=0=u_0(s_0)$, we achieve
 \bes\label{3.15}
 w(0,x)\geq u_0(x)\ \ \ {\rm on}\ \ [s_0-(2K)^{-1},\, s_0].
 \ees
Moreover, using the concavity of $w(0,x)$ and $w_x(0,s_0-K^{-1})=0$, it yields
 \bess
w(0,x)\geq w(0,s_0-(2M)^{-1})=\frac 34M\geq\|u_0\|_{C^1([0,s_0])}K^{-1}\geq  u_0(x)
 \eess
for $s_0-K^{-1}\leq x\leq s_0-(2K)^{-1}$, combining this with (\ref{3.15}) implies (\ref{3.13}).

Similarly, define
 \[\tilde K=\max \left\{\frac{4\|v_0\|_{C^1([0,s_0])}}{3M},\ \ \left(\frac{M^{q-1}}{2d_2}\right)^{1/2}\right\},\]
and
  $$z(t,x)=M[2\tilde K(s(t)-x)-\tilde K^2(s(t)-x)^2]$$
over the region $D^{\tilde K}$, we can prove
  \bes\label{3.16}
  v_x(t,s(t))\geq -2M\tilde K\ \ \ {\rm in}\ \ (0,T).
 \ees
It follows from (\ref{3.14}) and (\ref{3.16}) that
 $$s'(t)=-\mu(u_x+\rho v_x)(t,s(t))\leq C,$$
where positive constant $C$ depends on $M,\,\|u_0\|_{C^1([0,s_0])}$ and $\|v_0\|_{C^1([0,s_0])}$, but don't depend on $T$. This completes the proof.\ \ \ \ \fbox{}

\begin{lemma}\label{l3.3} \, Assume that $p,q\geq 1$. Let $T_{\max}$ and $(u,v,s)$ be obtained by Lemma $\ref{l3.1}$. If $T_{\max}<+\infty$, then $(\ref{1.7})$ holds.
\end{lemma}

{\bf Proof.}\, It is readily seen that if one component of $(u,v)$ blows up at time $T_{\max}$, so does the other one. Thus we suppose by contradiction that both $u$ and $v$
are bounded for $t\in [0,T_{\max})$ and $x\in[0,s(t)]$, namely, there exists a
positive constant $M$ such that
 \bes
 u(t,x)\leq M,\ \ v(t,x)\leq M,\ \ \forall \ (t,x)\in[0,T_{\max})\times[0,s(t)].
 \lbl{3.17}\ees
In terms of Lemma \ref{l3.2}, there is a positive constant $C$ independent
of $T_{\max}$ so that
  \bes
  0\leq s'(t)\leq C,\ \ \ s_0\leq s(t)\leq s_0+Ct\leq s_0+CT_{\max},\ \
  \forall  \ t\in[0,T_{\max}).\lbl{3.18}\ees

We shall prove that $(u,v,s)$ can be extended to $[0,\,T_{\max}+\tau]$ for some $\tau>0$ and get a contradiction with the definition of $T_{\max}$. To this aim, we first estimate $\|u,v\|_{C^{1+\frac\alpha 2,2+\alpha}([\varepsilon,T_{\max})\times[0,s(t)])}$ for any given $0<\varepsilon<T_{\max}$. Similar to the proof of Theorem \ref{t1.2}, under the transformation (\ref{2.2}), the relation (\ref{2.3}) holds and
  \bes\label{3.19}
\left\{\begin{array}{ll}
w_t-d_1f(t)w_{yy}-g(t,y)w_y=z^p, &0<t<T_{\max},\ 0<y<1,\\[1.2mm]
z_t-d_2f(t)z_{yy}-g(t,y)z_y=w^q, &0<t<T_{\max},\ 0<y<1,\\[1.2mm]
w_y(t,0)=z_y(t,0)=w(t,1)=z(t,1)=0, \ \ &0<t<T_{\max},\\[1.2mm]
w(0,y)=u_0(s_0y),\ z(0,y)=v_0(s_0y),&0\leq y\leq 1,
\end{array}\right.
 \ees
where $f(t)$ and $g(t,y)$ are as in the proof of Theorem \ref{t1.2}.
Thanks to (\ref{3.17}) and (\ref{3.18}), in view of the standard $L^p$ theory for
parabolic equations we find that $w,z\in W_k^{1,2}((0,T_{\max})\times(0,1))$, and there exists a positive constant $C_1$, which
depends only $T_{\max},\,M,\,C,\,\|u_0\|_{W_k^2((0,s_0))}$ and $\|v_0\|_{W_k^2((0,s_0))}$, such that $\|w,\,z\|_{W_k^{1,2}((0,T_{\max})\times(0,1))}\leq C_1$. Hence, by Sobolev's embedding theorem, $w,\,z,\,w_y,\,z_y\in C^{\frac\alpha 2,\,\alpha}([0,T_{\max})\times[0,1])$, and there exists a positive
constant $C_2$ depending only on $C_1$, $T_{\max}$ and $T_{\max}^{-1}$ such that
 \bes
 \|w,\,z,\,w_y,\,z_y\|_{C^{\frac\alpha 2,\,\alpha}([0,T_{\max})\times[0,1])}\leq C_2.
 \lbl{3.20}\ees
By use of (\ref{3.18}), (\ref{2.2}) and $u_x=s^{-1}(t)w_y,\,v_x=s^{-1}(t)z_y$, it follows from (\ref{3.20}) that
 \[\|u,\,v,\,u_x,\,v_x\|_{C^{\frac\alpha 2,\,\alpha}([0,T_{\max})\times[0,s(t)])}\leq C_3,\]
where $C_3$ depends only on $C_2$, $C$ and $T_{\max}$.
Therefore, by $s'(t)=-\mu\big[u_x(t,s(t))+\rho v_x(t,s(t))\big]$,
  \bes
  \|s'\|_{C^{\frac\alpha 2}([0,T_{\max}))}\leq\mu(1+\rho)C_3.\lbl{3.21}\ees

Fix $0<\varepsilon<T_{\max}$. Remember (\ref{3.18}), (\ref{3.20}) and (\ref{3.21}), we can apply the Schauder theory to problem (\ref{3.19}) in $[\varepsilon,T_{\max})\times[0,1]$, and obtain that $w,\,z\in C^{1+\frac\alpha 2,2+\alpha}\big([\varepsilon,T_{\max})\times[0,1]\big)$, and
 \bes
  \|w,\,z\|_{C^{1+\frac\alpha 2,2+\alpha}([\varepsilon,T_{\max})\times[0,1])} \leq C_4,\lbl{3.22}\ees
where $C_4$ is independent of $T_{\max}$.
Thanks to (\ref{3.18}), (\ref{2.2}) and (\ref{2.3}), it follows from (\ref{3.22}) that
 \bes\|u,\,v\|_{C^{1+\frac\alpha 2,2+\alpha}([\varepsilon,T_{\max})\times[0,s(t)])}\leq C_5
 \lbl{3.23}\ees
for some positive constant $C_5$ independent of $T_{\max}$.

Keeping in mind (\ref{3.23}) and following the proof of Lemma \ref{l3.1}, there exists a constant $\tau>0$ depending on $C_3$ and $C_5$ but independent of $T_{\max}$, such that the solution $(u,v,s)$ of (\ref{1.1}) with initial time $T_{\max}-\tau$ can be extended to the interval $[0,T_{\max}-\tau+2\tau]=[0,T_{\max}+\tau]$. The proof is finished.
\ \ \ \ \fbox{}

\section{Two comparison principles}\lbl{s.4}

In this section we present two comparison principles which play an important role in establishing the existence and uniqueness of maximal positive solution to (\ref{1.1}) when either $p<1$ or $q<1$.

\begin{lemma}\label{l4.1}\, Let $\underline a,\underline b$, $\bar a,\bar b$ be non-negative constants and $T\in(0,+\infty)$. Suppose that $\underline s,\,\bar s\in C^1([0,T])$ are positive functions in $[0,T]$, $\underline u,\,\underline v\in C(\overline Q_T)\cap C^{1,2}(Q_T)$ are positive functions in $Q_T$, and $\bar u,\,\bar v\in C(\overline Q^*_T)\cap C^{1,2}(Q^*_T)$ are positive in $Q^*_T$, where
 \bess
 Q_T&=&\{(t,x)\in\mathbb{R}^2:\ 0<t\leq T,\ 0<x<\underline s(t)\}, \\[1mm] Q^*_T&=&\{(t,x)\in\mathbb{R}^2:\ 0<t\leq T,\ 0<x<\bar s(t)\}.\eess
Assume further that $(\underline u,\,\underline v,\,\underline s)$ and $(\bar u,\,\bar v,\,\bar s)$ satisfy, in the classical sense,
 \bes
\left\{\begin{array}{ll}
\underline u_t-d_1\underline u_{xx}\leq (\underline v+\underline b)^p, &0<t\leq T,\ 0<x<\underline s(t),\\[1mm]
 \underline v_t-d_2 \underline v_{xx}\leq (\underline u+\underline a)^q, &0<t\leq T,\ 0<x<\underline s(t),\\[1mm]
 \underline s'(t)\leq-\mu (\underline u_x+\rho \underline v_x),&0<t\leq T,\ x=\underline s(t),\\[1mm]
 \underline u>0, \ \underline v>0, \ \underline u_x\geq 0,\ \underline v_x\geq0, \ \ &0<t\leq T, \ x=0,\\[1mm]
\underline u=\underline v=0,\ \ \ \ &0\leq t\leq T, \ x=\underline s(t)
\end{array}\right.\lbl{4.1}
 \ees
and
 \bes
\left\{\begin{array}{ll}
\bar u_t-d_1\bar u_{xx}\geq (\bar v+\bar b)^p, &0<t\leq T,\ 0<x<\bar s(t),\\[1mm]
 \bar v_t-d_2\bar v_{xx}\geq (\bar u+\bar a)^q, &0<t\leq T,\ 0<x<\bar s(t),\\[1mm]
 \bar s'(t)\geq-\mu (\bar u_x+\rho \bar v_x),&0<t\leq T,\ x=\bar s(t),\\[1mm]
\bar u>0, \ \bar v>0, \ \bar u_x\leq 0,\ \bar v_x\leq0, \ \ &0<t\leq T, \ x=0,\\[1mm]
\bar u=\bar v=0,\ \ \ \ &0\leq t\leq T, \ x=\bar s(t),
\end{array}\right.\lbl{4.2}
\ees
respectively. If $\underline a\leq\bar a$, $\underline b\leq\bar b$, $\bar u(0,\underline s(0))>0$, $\bar v(0,\underline s(0))>0$, $\underline s(0)<\bar s(0)$ and
 $$0<\underline u(0,x)\leq\bar u(0,x) ,\ 0<\underline v(0,x)\leq\bar v(0,x),\ \ \forall \ 0\leq x<\underline s(0),$$
then
 $$\big(\underline s(t),\,\underline u(t,x),\,\underline v(t,x)\big)<\big(\bar s(t),\,\bar u(t,x),\,\bar v(t,x)\big),\ \ \forall \ t\in(0,T],\ x\in [0,\underline s(t)].$$
\end{lemma}

{\bf Proof.} \, First of all we assert that $\underline s(t)<\bar s(t)$ for $0\leq t\leq T$. Obviously, this is true for small $t$. If this assertion does not hold, we can find a first $\tau<T$ so that $\underline s(t)<\bar s(t)$ for $t\in [0,\tau)$ and $\underline s(\tau)=\bar s(\tau)$. Thus
 \bes\label{4.3}
\underline s'(\tau)\geq \bar s'(\tau).
 \ees

Now we compare $(\underline u,\underline v)$ and $(\bar u,\bar v)$ over
  $$\Omega_{\tau}=\{(t,x)\in\mathbb{R}^2:\ 0<t<\tau,\ 0\leq x\leq\underline s(t)\}.$$
Note that $\bar u(t,\underline s(t))>0=\underline u(t,\underline s(t))$, $\bar v(t,\underline s(t))>0=\underline v(t,\underline s(t))$ for all $0\leq t<\tau$. For any given $0<\varepsilon\ll 1$, in terms of continuity, there exists a constant $0<\sigma_0\ll 1$ such that, for all $0<\sigma\leq\sigma_0$,
 \[\bar u(t,\underline s(t)-\sigma)>\underline u(t,\underline s(t)-\sigma)>0,\ \ \bar v(t,\underline s(t)-\sigma)>\underline v(t,\underline s(t)-\sigma)>0, \ \ \forall \ 0\leq t\leq\tau-\varepsilon.\]
Hence, $\bar u$, $\bar v$, $\underline u$ and $\underline v$ are positive in the domain
   \[\Omega^\varepsilon_\sigma=\{(t,x)\in\mathbb{R}^2:\ 0\leq t\leq\tau-\varepsilon,\ 0\leq x\leq \underline s(t)-\sigma\},\]
and $\bar u,\,\bar v,\,\underline u,\,\underline v\geq\delta$ in $\Omega^\varepsilon_\sigma$ for some constant $\delta>0$. Consequently, functions $(\underline v+\underline b)^p$, $(\underline u+\underline a)^q$, $(\bar v+\bar b)^p$ and $(\bar u+\bar a)^q$ are Lipschitz continuous when $(t,x)\in \Omega^\varepsilon_\sigma$. We can now apply the comparison principle to $(\bar u,\,\bar v)$ and $(\underline u,\,\underline v)$ in the domain $\Omega^\varepsilon_\sigma$, and conclude that $(\bar u,\,\bar v)\geq (\underline u,\,\underline v)$ in $\Omega^\varepsilon_\sigma$ since $\bar u(t,\underline s(t)-\sigma)>\underline u(t,\underline s(t)-\sigma)$, $\bar v(t,\underline s(t)-\sigma)>\underline v(t,\underline s(t)-\sigma)$. Letting $\sigma\to 0$, it is derived that $(\bar u,\,\bar v)\geq (\underline u,\,\underline v)$ in the domain
 \[\Omega^\varepsilon_0=\{(t,x)\in\mathbb{R}^2:\ 0\leq t\leq\tau-\varepsilon,\ 0\leq x\leq \underline s(t)\}.\]
By the arbitrariness of $\varepsilon>0$, we get
  \bes
  \big(\underline u(t,x),\,\underline v(t,x)\big)\leq\big(\bar u(t,x),\,\bar v(t,x)\big),\ \ \forall \ (t,x)\in\Omega_{\tau}.
  \lbl{4.4}\ees

Let $w=\bar u-\underline u$ and $z=\bar v-\underline v$. As $\underline a\leq\bar a$, $\underline b\leq\bar b$ and $p,q>0$, it follows from (\ref{4.1}) and (\ref{4.2}) that
 \bess
\left\{\begin{array}{ll}
w_t-d_1w_{xx}\geq (\bar v+\bar b)^p-(\underline v+\underline b)^p\geq 0, \ \ &0<t\leq\tau,\ 0<x<\underline s(t),\\[1mm]
 z_t-d_2 z_{xx}\geq (\bar u+\bar a)^q-(\underline u+\underline a)^q\geq 0, &0<t\leq\tau,\ 0<x<\underline s(t),\\[1mm]
 w=\bar u>0,\ z=\bar v>0, \ \ &0<t<\tau, \ x=\underline s(t),\\[1mm]
 w_x\leq 0,\ z_x\leq0, \ \ &0<t\leq\tau, \ x=0,\\[1mm]
w(0, x)\geq 0, \ \ z(0,x)\geq 0, \ &0\leq x\leq \underline s(0).
\end{array}\right.
 \eess
Therefore, $w>0$, $z>0$ in $\Omega_\tau$ by the strong maximum principle. Although we only know $\ud s\in C^1$, which indicates that the domain $\Omega_\tau$ may not have the {\it interior sphere} property at the right boundary $x=\ud s(t)$, in the same way as the proof of Theorem \ref{t1.3} for the case that $p<1$ or $q<1$, one can still find
$w_x(\tau,\underline s(\tau))<0$, $z_x(\tau,\underline s(\tau))<0$. As a conclusion,
$\bar u_x(\tau,\underline s(\tau))<\underline u_x(\tau,\underline s(\tau))$, $\bar v_x(\tau,\underline s(\tau))<\underline v_x(\tau,\underline s(\tau))$,
  $$\underline s'(\tau)\leq-\mu(\underline u_x+\rho \underline v_x)(\tau,\underline s(\tau))<-\mu(\bar u_x+\rho \bar v_x)(\tau,\underline s(\tau))\leq\bar s'(\tau).$$
This contradicts with (\ref{4.3}). So, $\underline s(t)<\bar s(t)$ for $0\leq t\leq T$.

The same argument as the proof of (\ref{4.4}) gives
$\underline u\leq\bar u$ and $\underline v\leq\bar v$ in $\Omega_{T}$. Hence, functions $w=\bar u-\underline u$ and $z=\bar v-\underline v$ satisfy
  \bess
\left\{\begin{array}{ll}
w_t-d_1w_{xx}\geq (\bar v+\bar b)^p-(\underline v+\underline b)^p\geq 0, \ \ &0<t\leq T,\ 0<x<\underline s(t),\\[1mm]
 z_t-d_2 z_{xx}\geq (\bar u+\bar a)^q-(\underline u+\underline a)^q\geq 0, &0<t\leq T,\ 0<x<\underline s(t),\\[1mm]
 w>0,\ z>0, \ \ &0<t\leq T, \ x=\underline s(t),\\[1mm]
 w_x\leq 0,\ z_x\leq0, \ \ &0<t\leq T, \ x=0,\\[1mm]
w(0, x)\geq 0, \ \ z(0,x)\geq 0, \ &0\leq x\leq \underline s(0).
\end{array}\right.
 \eess
By the strong maximum principle, $w>0$, $z>0$, as a result $\underline u<\bar u$, $\underline v<\bar v$ in $(0,T]\times[0,\underline s(t)]$.
The proof of Lemma \ref{l4.1} is complete. \ \ \ \ \fbox{}

\begin{lemma}\label{l4.2} Let $\underline a,\underline b$, $\bar a,\bar b$, $T$, $\underline s,\,\bar s$, $\underline u,\,\underline v$, and $\bar u,\,\bar v$ be as in Lemma $\ref{l4.1}$. Assume that $(\underline u,\,\underline v,\,\underline s)$ satisfies $(\ref{4.1})$ and $(\bar u,\,\bar v,\,\bar s)$ satisfies
 \bes
\left\{\begin{array}{ll}
\bar u_t-d_1\bar u_{xx}=(\bar v+\bar b)^p, &0<t\leq T,\ 0<x<\bar s(t),\\[1mm]
 \bar v_t-d_2\bar v_{xx}=(\bar u+\bar a)^q, &0<t\leq T,\ 0<x<\bar s(t),\\[1mm]
 \bar s'(t)=-\mu (\bar u_x+\rho \bar v_x),&0<t\leq T,\ x=\bar s(t),\\[1mm]
\bar u>0, \ \bar v>0, \ \bar u_x= 0,\ \bar v_x= 0, \ \ &0<t\leq T, \ x=0,\\[1mm]
\bar u=\bar v=0,\ \ \ \ &0\leq t\leq T, \ x=\bar s(t).
 \end{array}\right.\lbl{4.5}
 \ees
If $\underline a<\bar a$, $\underline b<\bar b$, $\underline s(0)\leq\bar s(0)$ and
 $$0<\underline u(0,x)\leq\bar u(0,x) ,\ 0<\underline v(0,x)\leq\bar v(0,x)\ \ \ {\rm for}\ \ 0\leq x<\underline s(0),$$
then
 \bes
 \big(\underline s(t),\,\underline u(t,x),\,\underline v(t,x)\big)<\big(\bar s(t),\,\bar u(t,x),\,\bar v(t,x)\big),\ \ \forall \ t\in(0,T],\ x\in [0,\underline s(t)].\quad
 \lbl{4.6}\ees
 \end{lemma}

{\bf Proof.}\, Since at least one reaction term is not locally Lipschitz continuous in (\ref{4.1}) if either (i) $\ud b=0$ and $p<1$, or (ii) $\ud a=0$ and $q<1$, the standard comparison principle cannot be used directly for these cases. To overcome this difficulty, we shall make approximations of initial functions and then apply Lemma \ref{l4.1} to get our desired conclusion.

Denote $\bar s_0=\bar s(0)$. For $0<\varepsilon\ll 1$, choose $u_\varepsilon$, $v_\varepsilon\in C^2([0,\bar s_0+\varepsilon])$ such that
 \bess
\left.\begin{array}{l}
 u_\varepsilon'(0)=v_\varepsilon'(0)=u_\varepsilon(\bar s_0+\varepsilon)=v_\varepsilon(\bar s_0+\varepsilon)=0,\\[1.5mm]
 u_\varepsilon(x), \ v_\varepsilon(x)>0 \ \ {\rm in} \ \ [0,\bar s_0+\varepsilon),\\[1.5mm]
u_\varepsilon(x)\geq\bar u(0,x),\ \ v_\varepsilon(x)\geq\bar v(0,x)\ \ {\rm on}\ \ [0,\bar s_0],\\[1.5mm]
(u_\varepsilon(x),v_\varepsilon(x)\big)\to(\bar u(0,x),\bar v(0,x)\big)\ \ \ {\rm in}\ \ W^2_k((0,\bar s_0))\ \ {\rm as}\ \ \varepsilon\to 0.
\end{array}\right.
 \eess
As $\bar a,\,\bar b>0$, similar to section \ref{s.3}, the problem
 \bess
\left\{\begin{array}{ll}
u_t-d_1u_{xx}=(v+\bar b)^p, \ \ &t>0,\ 0<x<s(t),\\[1mm]
v_t-d_2v_{xx}=(u+\bar a)^q, &t>0,\ 0<x<s(t),\\[1mm]
s'(t)=-\mu (u_x+\rho v_x),&t>0,\ x=s(t),\\[1mm]
u_x=v_x=0,&t>0, \ x=0,\\[1mm]
u=v=0,\ \ \ \ &t>0, \ x=s(t),\\[1mm]
u=u_\varepsilon,\ v=v_\varepsilon,\ \ &t=0, \ 0\leq x\leq\bar s_0+\varepsilon,\\[1mm]
s(0)=\bar s_0+\varepsilon
\end{array}\right.
 \eess
has a unique positive solution $(u_\varepsilon,\,v_\varepsilon,\,s_\varepsilon)$ defined in $[0,T_\varepsilon)$, here $T_\varepsilon$ is the maximum existence time of $(u_\varepsilon,\,v_\varepsilon,\,s_\varepsilon)$. Meanwhile, by Lemma \ref{l4.1},
 \bes
 \big(\underline s(t),\,\underline u(t,x),\,\underline v(t,x)\big)<\big(s_\varepsilon(t),\,u_\varepsilon(t,x),\, v_\varepsilon(t,x)\big)
 \lbl{4.7}\ees
for all $0<t<\min\{T_\varepsilon,T\}$ and $0\leq x\leq\underline s(t)$, and when $\varepsilon_1<\varepsilon_2$,
 \bes
 \big(s_{\varepsilon_1}(t),\,u_{\varepsilon_1}(t,x),\,v_{\varepsilon_1}(t,x)\big)< \big(s_{\varepsilon_2}(t),\,u_{\varepsilon_2}(t,x),\,v_{\varepsilon_2}(t,x)\big)\lbl{4.8}
 \ees
for all $0<t<\min\{T_{\varepsilon_1},\,T_{\varepsilon_2}\}$ and $0\leq x\leq s_{\varepsilon_1}(t)$.
Hence, $T_{\varepsilon_2}\leq T_{\varepsilon_1}$. There exist $\hat s(t)$, $\hat u(t,x)$, $\hat v(t,x)$ and $\hat T>0$ such that $\dd\lim_{\varepsilon\searrow 0}T_\varepsilon=\hat T$ and
 \bes
 \lim_{\varepsilon\searrow 0}\big(u_\varepsilon(t,x),\,v_\varepsilon(t,x),\,s_\varepsilon(t)\big)=\big(\hat u(t,x),\,\hat v(t,x),\,\hat s(t)\big)
 \lbl{4.9}\ees
for each $0\leq t<\hat T$ and $0\leq x\leq\hat s(t)$. Moreover,
 \[\big(\hat s(t),\,\hat u(t,x),\,\hat v(t,x)\big)<\big(s_\varepsilon(t),\,u_\varepsilon(t,x),\, v_\varepsilon(t,x)\big)\]
for all $0<t<T_\ep$, $0\leq x\leq \hat s(t)$ and $\varepsilon>0$ by (\ref{4.8}), and
 \bes
 \big(\underline s(t),\,\underline u(t,x),\,\underline v(t,x)\big)\leq(\hat s(t),\,\hat u(t,x),\,\hat v(t,x)\big)
\lbl{4.10}\ees
for all $0\leq t<\min\{\hat T,T\}$ and $0\leq x\leq\underline s(t)$ by (\ref{4.7}).

Since $\bar a,\,\bar b>0$, taking advantage of the $L^p$ estimate, Sobolev embedding theorem and interior Schauder estimate, it can be shown that, for any given $0<\sigma\ll 1$, there exists a positive constant $C_\sigma$ such that
 \bess
  &\|u_\ep\|_{C^{\frac{1+\alpha}2,1+\alpha}(Q_\varepsilon^\sigma)}
  +\|v_\ep\|_{C^{\frac{1+\alpha}2,1+\alpha}(Q_\varepsilon^\sigma)}
  +\|s_\ep\|_{C^{1+\frac\alpha 2}([0,T_\varepsilon-\sigma])}\leq C_\sigma,&\\[2mm]
  &\|u_\varepsilon\|_{C^{1+\frac\alpha 2,2+\alpha}\left(D_\varepsilon^\sigma\right)}
 +\|v_\varepsilon\|_{C^{1+\frac\alpha 2,2+\alpha}\left(D_\varepsilon^\sigma\right)}\leq C_\sigma,&\eess
where
 \bess
 & Q_\varepsilon^\sigma=\{(t,x)\in \mathbb{R}^2:\ t\in [0,T_\varepsilon-\sigma],\ x\in[0,s_\varepsilon(t)]\},&\\[1mm]
 &D_\varepsilon^\sigma=\{(t,x)\in \mathbb{R}^2:\ t\in [\sigma,T_\varepsilon-\sigma],\ x\in[0,s_\varepsilon(t)-\sigma]\}.&\eess
These estimates combined with (\ref{4.8}) and (\ref{4.9}) allow us to derive that $(u_\varepsilon,v_\varepsilon)\to(\hat u, \hat v)$ in $[C^{1,2}_{\rm loc}(D_\sigma)\cap C^{0,1}(\overline D_\sigma)]^2$ and $s_\varepsilon\to \hat s$ in $C^1([0,\hat T-\sigma])$ as $\varepsilon\to 0$, where
 \bess
 D_\sigma=\{(t,x)\in \mathbb{R}^2:\ t\in (0,\hat T-\sigma),\ x\in(0,\hat s(t))\}.\eess
Consequently, by the arbitrariness of $\sigma>0$, we observe that $(\hat u, \hat v, \hat s)$ satisfies
  \bes
\left\{\begin{array}{ll}
\hat u_t-d_1 \hat u_{xx}=(\hat v+\bar b)^p, &0<t<\hat T,\ 0<x<\hat s(t),\\[1mm]
 \hat v_t-d_2 \hat v_{xx}=(\hat u+\bar a)^q, &0<t<\hat T,\ 0<x<\hat s(t),\\[1mm]
 \hat s'(t)=-\mu (\hat u_x+\rho \hat v_x),&0<t<\hat T,\ x=\hat s(t),\\[1mm]
 \hat u_x= 0,\ \hat v_x= 0, \ \ &0<t<\hat T, \ x=0,\\[1mm]
 \hat u(t, \hat s(t))=\hat v(t, \hat s(t))=0,\ \ \ \ &0\leq t<\hat T,\\[1mm]
 \hat u(0,x)=\bar u(0,x), \ \hat v(0,x)=\bar v(0,x), \ &0\leq x\leq\bar s(0),\\[1mm]
 \hat s(0)=\bar s(0).
 \end{array}\right.\lbl{4.11}
 \ees
Since $\bar a,\,\bar b>0$, the uniqueness result is true for problem (\ref{4.11}). Obviously, $(\bar u,\,\bar v,\,\bar s)$ satisfies (\ref{4.11}). And so, $(\bar u,\,\bar v,\,\bar s)\equiv (\hat u,\,\hat v,\,\hat s)$ for $0\leq t\leq\min\{\hat T,T\}$. Taking into account (\ref{4.10}), one has
 \bes
 \big(\underline s(t),\,\underline u(t,x),\,\underline v(t,x)\big)\leq\big(\bar s(t),\,\bar u(t,x),\,\bar v(t,x)\big),
 \ \ \forall \ 0\leq t<\min\{\hat T,T\}, \ 0\leq x\leq\underline s(t).
 \lbl{4.12}\ees

Denote $T^*=\min\{\hat T,T\}$. We claim that
 \bes
 {\rm For \ any \ given} \ 0<\rho<T^*,\ {\rm there \ must \ be \ a } \ 0<\tau<\rho \ {\rm such \ that} \ \underline s(\tau)<\bar s(\tau).
 \lbl{4.13}\ees
If this is not true, then $\underline s(t)=\bar s(t)$ for all $0\leq t<\rho$. Let $w=\bar u-\underline u$ and $z=\bar v-\underline v$. As $\underline a<\bar a$, $\underline b<\bar b$ and $p,q>0$,  utilizing (\ref{4.1}), (\ref{4.5}) and (\ref{4.12}) one can derive that $w\geq 0$, $z\geq 0$ and satisfy
 \bess
\left\{\begin{array}{ll}
w_t-d_1w_{xx}\geq (\bar v+\bar b)^p-(\underline v+\underline b)^p>0, \ \ &0<t<\rho,\ 0<x<\underline s(t),\\[1mm]
 z_t-d_2 z_{xx}\geq (\bar u+\bar a)^q-(\underline u+\underline a)^q>0, &0<t<\rho,\ 0<x<\underline s(t),\\[1mm]
 w=z=0, \ \ &0<t<\rho, \ x=\underline s(t),\\[1mm]
 w_x\leq 0,\ z_x\leq0, \ \ &0<t<\rho, \ x=0,\\[1mm]
w(0, x)\geq 0, \ \ z(0,x)\geq 0, \ &0\leq x\leq \underline s(0).
\end{array}\right.
 \eess
By use of the strong maximum principle, $w(t,x)>0$, $z(t,x)>0$ for all $0<t<\rho$ and $0\leq x<\underline s(t)$. Similar to the argument in the proof of Lemma \ref{l4.1}, we acquire $w_x(t,\underline s(t))<0$, $z_x(t,\underline s(t))<0$, and then get $\bar u_x(t,\underline s(t))<\underline u_x(t,\underline s(t))$, $\bar v_x(t,\underline s(t))<\underline v_x(t,\underline s(t))$ for all $0<t<\rho$.  Thus,
 $$\underline s'(t)\leq-\mu(\underline u_x+\rho \underline v_x)(t,\underline s(t))<-\mu(\bar u_x+\rho \bar v_x)(t,\underline s(t))=\bar s'(t), \ \ \forall \ 0<t<\rho.$$
Integrating the above inequality over $[0,t]$, it yields
 \[\underline s(t)-\underline s(0)<\bar s(t)-\bar s(0) , \ \ \forall \ 0<t<\rho,\]
which implies $\underline s(t)<\bar s(t)$ for all $0<t<\rho$. This is a contradiction, and hence our claim (\ref{4.13}) is true.

Choose $t_n\searrow 0$ such that $\underline s(t_n)<\bar s(t_n)$. Then $(\underline u,\,\underline v,\,\underline s)$ satisfies $(\ref{4.1})$ and $(\bar u,\,\bar v,\,\bar s)$ satisfies (\ref{4.5}) where the interval $[0,T]$ is replaced by $[t_n,T]$. Remember $\underline a<\bar a$, $\underline b<\bar b$, $\underline s(t_n)<\bar s(t_n)$ and
 $$0<\underline u(t_n,x)\leq\bar u(t_n,x) ,\ 0<\underline v(t_n,x)\leq\bar v(t_n,x)\ \ {\rm for}\ \
 0\leq x<\underline s(t_n).$$
By virtue of Lemma \ref{l4.1}, it follows that
 $$\big(\underline s(t),\,\underline u(t,x),\,\underline v(t,x)\big)<\big(\bar s(t),\,\bar u(t,x),\,\bar v(t,x)\big),\ \ \forall \ t\in(t_n,T],\ x\in [0,\underline s(t)].$$
Letting $t_n\searrow 0$, the conclusion (\ref{4.6}) is obtained.
\ \ \ \ \fbox{}

\begin{remark}\lbl{r4.1} \, In Lemma $\ref{l4.2}$, if the conditions $\underline a<\bar a$ and $\underline b<\bar b$ are replaced by $0<\underline a\leq\bar a$, $0<\underline b\leq\bar b$ and either $\underline a<\bar a$ or $\underline b<\bar b$, then the conclusion is remains true.
\end{remark}

\section{Proof of Theorem \ref{t1.1}: the case with either $p<1$ or $q<1$}\lbl{s.5}

In this section we shall prove the following theorem.

\begin{theorem}\lbl{t5.1} \, Let $p<1$ or $q<1$. Then there exist a maximum existence time $T_{\max}$ and a unique maximal positive solution $(\bar u,\,\bar v,\,\bar s)$ of $(\ref{1.1})$ defined in $[0,T_{\max})$, such that either $T_{\max}=+\infty$, or $T_{\max}<+\infty$ and
 \bes\label{5.1}
  \limsup_{T\nearrow T_{\max}}\|\bar u\|_{L^{\infty}([0,T]\times[0,\bar s(t)])}=+\infty,\ \ \  \limsup_{T\nearrow T_{\max}}\|\bar v\|_{L^{\infty}([0,T]\times[0,\bar s(t)])}=+\infty.
 \ees
\end{theorem}

Let us point out that the existence in time of solution cannot be obtained directly by means of the contraction mapping theorem for our present situation, since at least one reaction term is not locally Lipschitz continuous in the unknown. To overcome this difficulty, we approximate the resource terms to get the local existence result.

{\bf Proof of Theorem \ref{t5.1}.}\, {\it Step 1}\, The construction of $(\bar u,\,\bar v,\,\bar s)$.

For $n\geq1$, we first consider the approximating problem
 \bes\label{5.2}
\left\{\begin{array}{ll}
u_t-d_1u_{xx}=(v+\frac 1n)^p, &t>0,\ 0<x<s(t),\\[1mm]
v_t-d_2v_{xx}=(u+\frac 1n)^q, &t>0,\ 0<x<s(t),\\[1mm]
s'(t)=-\mu (u_x+\rho v_x),&t>0,\ x=s(t),\\[1mm]
u_x=v_x=0,&t>0, \ x=0,\\[1mm]
u=v=0,\ \ \ \ &t>0, \ x=s(t),\\[1mm]
u(0,x)=u_0(x),\ v(0,x)=v_0(x),\ \ &0\leq x\leq s_0,\\[1mm]
s(0)=s_0.
\end{array}\right.
 \ees
Since functions $(v+\frac 1n)^p$ and $(u+\frac 1n)^q$ are Lipschitz continuous for $u,v\geq 0$, same as the arguments for the case $p\geq1,q\geq1$, it can be shown that there exist a maximum existence time $T_n>0$ and a unique positive solution $(u_n,v_n,s_n)$ of problem (\ref{5.2}) defined in $[0,T_n)$, such that either $T_n=+\ty$, or $T_n<+\ty$ and
 \bess
 \limsup_{T\nearrow T_n}\|u_n(t,x)\|_{L^{\infty}([0,T]\times[0,s_n(t)])}=+\infty,\ \ \  \limsup_{T\nearrow T_n}\|v(t,x)\|_{L^{\infty}([0,T]\times[0,s_n(t)])}=+\infty.
 \eess
Moreover, Lemma \ref{l4.2} gives
 \[\big(s_{n+1}(t),\,u_{n+1}(t,x),\,v_{n+1}(t,x)\big)<\big(s_n(t),\,u_n(t,x),\,v_n(t,x)\big)\]
for all $0<t<\min\{T_n,\,T_{n+1}\}$ and $x\in [0, s_{n+1}(t)]$, this implies $T_n=\min\{T_n,\,T_{n+1}\}$.

Similar to the arguments in the proof of Lemma \ref{l4.2}, there exists $\bar T\leq+\ty$ such that $T_n\to\bar T$, and there exist  $\bar s\in C^1([0,\bar T))$, $\bar u,\,\bar v\in C^{1,2}((0,\bar T)\times[0,\bar s(t))\cap C^{0,1}([0,\bar T)\times[0,\bar s(t))]$, such that, for any given $0<\sigma\ll 1$, $s_n\to\bar s$ in $C^1([0,\bar T-\sigma])$, $(u_n,v_n)\to (\bar u,\,\bar v)$ in $[C^{1,2}_{\rm loc}(Q_\sigma)\cap C^{0,1}(\overline Q_\sigma)]^2$, where
  \[Q_\sigma=\{(t,x)\in \mathbb{R}^2:\ t\in(0,\bar T-\sigma),\ x\in(0,\bar s(t))\}.\]
This asserts that $(\bar u,\,\bar v,\,\bar s)$ satisfies (\ref{1.1}) in $(0,\bar T)$ since  $\sigma$ is arbitrary. Moreover, for all $n$, we gain
 \bes
   (\bar s(t),\,\bar u(t,x),\,\bar v(t,x)\big)<(s_n(t),\,u_n(t,x),\,v_n(t,x)\big),\ \ \forall \ 0<t<T_n,\ x\in [0,\,\bar s(t)].
 \lbl{5.3}\ees

{\it Step 2}\, It will be proved that $(\bar u,\,\bar v,\,\bar s)$ is the maximal positive solution of $(\ref{1.1})$ defined in $[0,\bar T)$.

Obviously, $\bar u(t,x),\,\bar v(t,x)\geq 0$. Thanks to (\ref{1.1}), it yields
 \bess
\left\{\begin{array}{ll}
\bar u_t-d_1\bar u_{xx}=\bar v^p\geq 0, &0<t<\bar T,\ 0<x<\bar s(t),\\[1mm]
 \bar v_t-d_2\bar v_{xx}=\bar u^q\geq 0, &0<t<\bar T,\ 0<x<\bar s(t),\\[1mm]
 \bar u_x= 0,\ \bar v_x= 0, \ \ &0<t<\bar T, \ x=0,\\[1mm]
\bar u=\bar v=0,\ \ \ \ &0\leq t<\bar T, \ x=\bar s(t),\\[1mm]
\bar u(0,x)=u_0(x),\ \bar v(0,x)=v_0(x),\ \ &0\leq x\leq s_0,\\[1mm]
\bar s(0)=s_0.
 \end{array}\right.
 \eess
Since $u_0,\,v_0>0$ in $[0,s_0)$, it follows that
$\bar u,\,\bar v>0$ in $[0,\bar T)\times[0,\bar s(t))$ by the maximum principle. This suggests that  $(\bar u,\,\bar v,\,\bar s)$ is a positive solution of (\ref{1.1}) defined in $[0,\bar T)$.

Suppose that $(u,v,s)$ is a positive solution of $(\ref{1.1})$ defined in $[0,T]$ for some $0<T<\bar T$. Then there exists $n_0\gg1$ such that $T<T_n$ for all $n\geq n_0$. By use of Lemma \ref{l4.2} we have that, for all $n\geq n_0$,
 \[\big(s(t),\,u(t,x),\,v(t,x)\big)<\big(s_n(t),\,u_n(t,x),\,v_n(t,x)\big),\ \ \forall \ 0<t<T,\ x\in [0, s(t)].\]
Letting $n\to+\ty$, it yields
   \[(s(t),\,u(t,x),\,v(t,x)\big)\leq\big(\bar s(t),\,\bar u(t,x),\,\bar v(t,x)\big),\ \ \ \forall \ 0<t<T,\ x\in [0, s(t)].\]

{\it Step 3}\, In this step we shall prove that if $\bar T<+\ty$ and $\bar u,\,\bar v$ are bounded in $[0,\bar T)\times[0,\bar s(t)]$, then $(\bar u,\,\bar v,\,\bar s)$ can be extended to $[0,\bar T+\tau)$ for some $\tau>0$ such that  $(\bar u,\,\bar v,\,\bar s)$ is the maximal positive solution of $(\ref{1.1})$ defined in $[0,\bar T+\tau)$.

Choose positive constant $K$ such that
 $$\bar u(t,x)\leq K,\ \ \bar v(t,x)\leq K,\ \ \forall \ (t,x)\in[0,\bar T)\times[0,\bar s(t)].$$
Thanks to Lemma \ref{l3.2}, there exists a positive constant $C$ independent of $\bar T$ such that $0<\bar s'(t)\leq C$ for $t\in(0,\bar T)$. In the same way as the proof of Lemma \ref{l3.3}, it can be shown that, for any $0<\varepsilon\ll 1$,
 \[\bar s\in C^{1+\frac{1+\beta} 2}([\varepsilon,\bar T]), \ \ \bar u,\,\bar v\in C^{1+\frac \beta 2,\,2+\beta}([\varepsilon,\bar T)\times[0,\bar s(t)]),\]
and
 \[\|\bar s'\|_{C^{\frac{1+\beta} 2}([\varepsilon,\bar T])}\leq C_6,\ \
 \|\bar u,\,\bar v\|_{C^{1+\frac\beta 2,\,2+\beta}([\varepsilon,\bar T)\times[0,\bar s(t)])}
 \leq C_6\]
for some positive constant $C_6$ independent of $\bar T$, where $\beta=\min\{p,\,q\}$.
Following the proof of step 1, we assert that there exists $\tau>0$ such that, for any fixed $T_0: \varepsilon<T_0<\bar T$, the problem
 \bes
\left\{\begin{array}{ll}
u_t-d_1u_{xx}=v^p, &T_0<t<T_0+2\tau,\ 0<x<s(t),\\[1mm]
v_t-d_2v_{xx}=u^q, &T_0<t<T_0+2\tau,\ 0<x<s(t),\\[1mm]
s'(t)=-\mu (u_x+\rho v_x),\ \ &T_0<t<T_0+2\tau,\ x=s(t),\\[1mm]
u_x=v_x=0,&T_0<t<T_0+2\tau, \ x=0,\\[1mm]
u=v=0,\ \ \ \ &T_0<t<T_0+2\tau,\ x=s(t),\\[1mm]
u=\bar u,\ v=\bar v,\ \ &t=T_0, \ 0\leq x\leq\bar s(T_0),\\[1mm]
s(T_0)=\bar s(T_0)
\end{array}\right.\lbl{5.4}
 \ees
has a maximal positive solution defined in $[T_0,T_0+2\tau]$. Take $T_0=\bar T-\tau$, and denote the corresponding maximal positive solution of (\ref{5.4}) by $(u^*,\,v^*,\,s^*)$, which is defined in $[\bar T-\tau,\,\bar T+\tau]$.
It is obvious that $(\bar u,\,\bar v,\,\bar s)$ satisfies the equations of (\ref{5.4}) in $[T_0,\,\bar T)$. Because $(u^*,\,v^*,\,s^*)$ is the maximal positive solution of (\ref{5.4}) defined in $[T_0,\,\bar T+\tau]$, we find that
 \bes
 \big(\bar s(t),\,\bar u(t,x),\,\bar v(t,x)\big)\leq\big(s^*(t),\,u^*(t,x),\,v^*(t,x)\big),\ \ \forall \ T_0\leq t<\bar T,\ x\in [0,\,\bar s(t)].
 \lbl{5.5}\ees

Let $(u_n,v_n,s_n)$ be the unique positive solution of (\ref{5.2}) defined in $[0,T_n)$. Since $T_0<\bar T$, there exists $n_1\gg 1$ such that $T_n>T_0+\sigma$ for some $\sigma>0$ and all $n\geq n_1$. It is obvious that $(u_n,v_n,s_n)$ satisfies (for $n\geq n_1$)
 \bes
\left\{\begin{array}{ll}
u_t-d_1u_{xx}=(v+\frac 1n)^p, &T_0<t<T_0+\sigma,\ 0<x<s(t),\\[2mm]
v_t-d_2v_{xx}=(u+\frac 1n)^q, &T_0<t<T_0+\sigma,\ 0<x<s(t),\\[2mm]
s'(t)=-\mu (u_x+\rho v_x),&T_0<t<T_0+\sigma,\ x=s(t),\\[1.5mm]
u_x(t,0)=v_x(t,0)=0,&T_0<t<T_0+\sigma,\\[1.5mm]
u(t,s(t))=v(t,s(t))=0,\ \ \ \ &T_0<t<T_0+\sigma,\\[1.5mm]
u(T_0,x)=u_n(T_0,x),\ v(T_0,x)=v_n(T_0,x),\ \ &0\leq x\leq s_n(T_0),\\[1.5mm]
s(T_0)=s_n(T_0)
\end{array}\right.\lbl{5.6}
 \ees
Taking advantage of (\ref{5.3}), we have $\bar s(T_0)<s_n(T_0)$ and
 \[\big(\bar u(T_0,x),\,\bar v(T_0,x)\big)\leq\big(u_n(T_0,x),\,v_n(T_0,x)\big),\ \ \forall \ x\in [0,\bar s(T_0)].\]
Applying Lemma \ref{l4.1} to (\ref{5.4}) and (\ref{5.6}), it can be deduced that
 \[\big(s^*(t),\,u^*(t,x),\,v^*(t,x)\big)<\big(s_n(t),\,u_n(t,x),\,v_n(t,x)\big),\ \ \forall \ T_0\leq t<T_0+\sigma,\ x\in [0, s^*(t)].\]
Letting $n\to+\ty$, it yields
  \[\big(s^*(t),\, u^*(t,x),\,v^*(t,x)\big)\leq\big(\bar s(t),\,\bar u(t,x),\,\bar v(t,x)\big),\ \ \forall \ T_0\leq t<T_0+\sigma,\ x\in [0, s^*(t)].\]
This combined with (\ref{5.5}) allows us to derive that
 \[\big(\bar s(t),\,\bar u(t,x),\,\bar v(t,x)\big)=\big(s^*(t),\,u^*(t,x),\,v^*(t,x)\big),\ \ \forall \ T_0\leq t<T_0+\sigma,\ x\in [0, s^*(t)].\]
Define
 \bess
 &\hat s(t)=\bar s(t), \ \ \hat u(t,x)=\bar u(t,x), \ \ \hat v(t,x)=\bar v(t,x) \ \ \
 {\rm for} \ \ 0\leq t<T_0+\sigma,\ x\in [0, s^*(t)],&\\[1mm]
 &\hat s(t)=s^*(t), \ \ \hat u(t,x)=u^*(t,x), \ \ \hat v(t,x)=v^*(t,x) \ \ \
 {\rm for} \ \ T_0\leq t<T_0+2\tau,\ x\in [0, s^*(t)],&
 \eess
then $(\hat u,\,\hat v,\,\hat s)$ is a positive solution of (\ref{1.1}) defined in $[0,\bar T+\tau]$ since $T_0+2\tau=\bar T+\tau$. Same as the argument of step 2, we can still prove that $(\hat u,\,\hat v,\,\hat s)$ is the maximal positive solution of (\ref{1.1}) defined in $[0,\bar T+\tau]$.

{\it Step 4} \, Finally, making the extension of $(\hat u,\,\hat v,\,\hat s)$ step by step to a large existence interval as in Step 3, one can get a $T_{\max}$ and the unique maximal positive solution $(u,v,s)$ of (\ref{1.1}) defined in $[0,T_{\max})$ such that either $T_{\max}=+\infty$, or $T_{\max}<+\infty$ and (\ref{5.1}) holds. The proof is complete. \ \ \ \ \fbox{}

\section{Global existence, finite time blow-up and long time behavior: proofs of Theorems \ref{t1.4} and \ref{t1.5}}\label{s.6}

In this section, we study global existence and finite time blow-up of positive solution to (\ref{1.1}), and get long time behavior of bounded global solution. That is, to prove Theorems \ref{t1.4} and \ref{t1.5}. For the convenience to readers, we repeat them here.

\begin{theorem}\label{t6.1} \, Let $s_0$, $\mu$ and $\rho$ be fixed, $(u,v,s)$ and $T_{\max}$ be obtained in Theorem $\ref{t1.1}$.

{\rm(i)} Assume that $pq>1$. Then $T_{\max}=+\ty$, i.e., $(u,v)$ exists globally in time provided initial functions $u_0(x)$ and $v_0(x)$ are suitably small; while $T_{\max}<+\ty$, i.e., $(u,v)$ will blow up in finite time provided initial functions $u_0(x)$ and $v_0(x)$ are large enough.

{\rm(ii)} Assume that $0<pq\leq1$, then $T_{\max}=+\ty$.
\end{theorem}

{\bf Proof.} \, (i) We shall use the argument from Ricci and Tarzia (\cite{RT}) to construct a suitable upper solution and use Lemma \ref{l4.1} to derive that $T_{\max}=+\ty$ if the initial functions $u_0(x)$ and $v_0(x)$ are suitably small. Assume for definiteness that $0<p\leq1<q$. Since $pq>1$, there exist $\varepsilon_1,\varepsilon_2$ satisfying
 \[0<\varepsilon_1,\varepsilon_2<\frac d{8\mu(1+\rho)}\]
such that
 \bes\label{6.1}
{\varepsilon_1d}-16s_0^2\varepsilon_2^p\geq0,\ \ \
{\varepsilon_2d}-16s_0^2\varepsilon_1^q\geq0,
 \ees
where $d=\min\{d_1,d_2\}$. For such fixed positive constants $\varepsilon_1,\varepsilon_2$, define
  $$\bar s(t)=2s_0(2-{\rm e}^{-a t}),\ \ t\geq0;\ \ \ \ w(x)=1-y^2,\ \ 0\leq y\leq1,$$
and
   $$\bar u(t,x)=\varepsilon_1{\rm e}^{-b t}w\left(\frac x{\bar s(t)}\right),\ \ \bar v(t,x)=\varepsilon_2{\rm e}^{-\gamma t}w\left(\frac x{\bar s(t)}\right),\ \ t\geq0 ,\ \ 0\leq x\leq\bar s(t),$$
where $a,\,b$ and $\gamma$ are real parameters to be chosen later.

It will be shown that if $(u_0,\,v_0)$ satisfies
  \bes\label{6.2}
 \|u_0\|_{L^\infty}\leq\varepsilon_1/2,\ \ \ \|v_0\|_{L^\infty}\leq\varepsilon_2/2,
 \ees
then
 \bes\label{6.3}
   \big(s(t),\,u(t,x),\,v(t,x)\big)<\big(\bar s(t),\,\bar u(t,x),\,\bar v(t,x)\big),\ \ \forall \ 0\leq t<T_{\max},\ \ x\in [0,s(t)).
 \ees
And hence, $T_{\max}=+\ty$ by conclusion (\ref{1.7}).

We shall prove the estimate (\ref{6.3}) by use of Lemma \ref{l4.1}. Thanks to $0\leq w\leq 1$, $w'(y)=-2y\leq 0$, $w''(y)=-2$ and $\bar s'(t)>0$, elementary calculations give that, for all $t>0$ and $0<x<\bar s(t)$,
\bess
  \bar u_t-d_1\bar u_{xx}-\bar v^p &=&\dd\varepsilon_1{\rm e}^{-b t}\left(-b w-w'x\bar s'\bar s^{-2}-d_1w''\bar s^{-2}\right)-\varepsilon_2^pe^{-\gamma pt}w^p \\[1mm]
  &\geq& \varepsilon_1{\rm e}^{-b t}\dd\left(-b-w'x\bar s'\bar s^{-2}-d_1w''\bar s^{-2}\right)-\varepsilon_2^pe^{-\gamma pt}w^p \\[1mm]
   &\geq& \varepsilon_1{\rm e}^{-b t}\dd\left(-b+{2d\bar s^{-2}}\right)-\varepsilon_2^pe^{-\gamma pt}\\[1mm]
   &\geq& \varepsilon_1{\rm e}^{-b t}\dd\left(-b+{8^{-1}d\bar s_0^{-2}}\right)-\varepsilon_2^pe^{-\gamma pt}.
  \eess
Choose $b=\gamma p={d/16s_0^2}$. It follows that, for all $t>0$ and $0<x<\bar s(t)$,
 \bess\label{d3}
\bar u_t-d_1\bar u_{xx}-\bar v^p\geq {\rm e}^{-b t}\left({\varepsilon_1d(16s_0^2)^{-1}}-\varepsilon_2^p\right)\geq 0
 \eess
by the first inequality of (\ref{6.1}). Similarly, for all $t>0$ and $0<x<\bar s(t)$,
 \bess\label{d4}
\bar v_t-d_2\bar v_{xx}-\bar u^q\geq {\rm e}^{-\gamma t}\left({\varepsilon_2d(16s_0^2)^{-1}}-\varepsilon_1^q\right) \geq 0
 \eess
holds. Set $\varepsilon=\max\{\varepsilon_1,\varepsilon_2\}$ and choose $a=b$, then
 \begin{eqnarray*}
\bar s'(t)+\mu(\bar u_x+\rho\bar v_x)(t,\bar s(t))&=& 2a s_0{\rm e}^{-a t}-2\dd\mu\left(\varepsilon_1{\rm e}^{-b t}+\rho\varepsilon_2{\rm e}^{-\gamma t}\right)\bar s^{-1}(t)\nonumber\\
  &\geq& \dd 8^{-1}{ {\rm e}^{-a t}s_0^{-1}}\left[d-8{\mu(1+\rho)\varepsilon}\right]>0.
 \end{eqnarray*}
On the other hand, it is easy to show that
 \[\bar u_x(t,0)=\bar v_x(t,0)=0,\ \
 \bar u(t,\bar s(t))=\bar v(t,\bar s(t))=0,\ \ \ \forall\ t>0,\]
and by (\ref{6.2})
 \[u_0(x)<\bar u(0,x),\ \ v_0(x)<\bar v(0,x),\ \ \forall \ 0\leq x\leq s_0.\]
Since $\bar s(0)=2s_0>s_0$, thanks to Lemma \ref{l4.1}, we can derive (\ref{6.3}).

Now we show that $T_{\max}<+\ty$ provided that the initial functions $u_0(x)$ and $v_0(x)$ are large enough. Let $s^*=s_0/2$, $w_0(x)=u_0(x)/2$ and $z_0(x)=v_0(x)/2$. Since $pq>1$, it is well known that (cf. \cite{EH3}) the solution $(w,\,z)$ of
 \bes\label{6.4}
\left\{\begin{array}{ll}
w_t-d_1w_{xx}=z^p, &t>0,\ 0<x<s^*,\\[1mm]
z_t-d_2z_{xx}=w^q, &t>0,\ 0<x<s^*,\\[1mm]
w_x(t,0)=z_x(t,0)=0,&t>0,\\[1mm]
w(t,s^*)=w(t,s^*)=0,\ \ \ \ &t>0,\\[1mm]
w(0,x)=w_0(x),\ z(0,x)=z_0(x),\ \ &0\leq x\leq s^*\\[1mm]
\end{array}\right.
 \ees
will blow up in finite time provided that $u_0(x)$ and $v_0(x)$ are large enough. Thanks to $s(t)>s^*$ for all $0\leq t<T_{\max}$, we can apply the comparison principle (Lemma 2.2 of \cite{EH3}) to (\ref{1.1}) and (\ref{6.4}) and get that $u(t,x)\geq w(t,x)$, $v(t,x)\geq z(t,x)$. Thus, $(u,v)$ blows up in finite time, and in turn $T_{\max}<+\ty$.

(ii) We make the zero extension of $u_0(x)$ and $v_0(x)$ to $(s_0, +\ty)$, and consider
\bes\label{6.5}
\left\{\begin{array}{ll}
w_t-d_1w_{xx}=z^p,\ \ &t>0,\ \ x>0,\\[1mm]
z_t-d_2z_{xx}=w^q,\ \ &t>0,\ \ x>0,\\[1mm]
w_x(t,0)=z_x(t,0)=0,\ \ &t>0,\\[1mm]
w=u_0+1,\ \ z=v_0+1,\ \ &t=0, \ x\geq 0.
\end{array}\right.
  \ees
It is well known (see, Theorem 1 in \cite{EH2}) that solution $(w,z)$ to  problem (\ref{6.5}) exists globally in time.  By the comparison principle (Lemma 2.2 of \cite{EH3}) one has that $u(t,x)\leq w(t,x)$, $v(t,x)\leq z(t,x)$. Therefore, in view of (\ref{1.7}), we deduce that $T_{\max}=+\ty$.

This completes the proof of Theorem \ref{t6.1}
\ \ \ \ \fbox{}

In what follows, we demonstrate the asymptotic behavior of bounded global solution. To this aim, we first give a lemma which can be proved by the similar way to that of Proposition 3.1 in \cite{WMX} and the details will be omitted.

\begin{lemma}\lbl{l6.1} Let $d$, $\beta$ and $s_0$ be positive constants and $C\in\mathbb{R}$. Assume that functions $s(t)$ and $w(t,x)$  satisfy $s(t)>0$,  $w(t,x)>0$ for all $0\leq t<+\infty$ and $0<x<s(t)$. We further suppose that
  $$ \lim_{t\to+\infty} s(t)=s_\infty<+\infty, \ \ \lim_{t\to+\infty} s'(t)=0,$$
and
 $$ \|w(t,\cdot)\|_{C^1[0,\,s(t)]}\leq M, \ \ \forall \ t>1$$
for some constant $M>0$. If $(w,s)$ satisfies
  \bess\left\{\begin{array}{lll}
 w_t-dw_{xx}\geq Cw, &t>0,\ \ 0<x<s(t),\\[.5mm]
 w_x=0,\ {\rm or} \ w=0, \ &t>0, \ \ x=0,\\[.5mm]
 w=0,\ \ s'(t)\geq-\beta w_x, \ &t>0,\ \ x=s(t),\\[.5mm]
 w(0,x)=w_0(x), \ &x\in [0,s_0],\\[.5mm]
 s(0)=s_0
 \end{array}\right.\eess
in the classical sense, then $\dd\lim_{t\to+\infty}\,\max_{0\leq x\leq s(t)}w(t,x)=0$.
\end{lemma}

\begin{theorem}\lbl{t6.2}\, Let $s_0$, $\mu$ and $\rho$ be fixed, $(u,v,s)$ and $T_{\max}$ be obtained in Theorem $\ref{t1.1}$. If $T_{\max}=+\ty$, $s_{\infty}:=\lim_{t\to+\infty}s(t)<+\infty$, $u$ and $v$ are bounded, then
 \bes
 \lim_{t\to+\infty}\,\max_{0\leq x\leq s(t)}u(t,x)=\lim_{t\to+\infty}\,\max_{0\leq x\leq s(t)}v(t,x)=0.
 \lbl{6.6}\ees
\end{theorem}

{\bf Proof.} \, Under the transformation (\ref{2.2}), $(w,z)$ satisfies
  \bes\label{6.7}
\left\{\begin{array}{ll}
w_t-d_1f(t)w_{yy}-g(t,y)w_y=z^p, &t>0,\ 0<y<1,\\[1mm]
z_t-d_2f(t)z_{yy}-g(t,y)z_y=w^q, &t>0,\ 0<y<1,\\[1mm]
w_y(t,0)=z_y(t,0)=w(t,1)=z(t,1)=0, \ \ &t>0,\\[1mm]
w(0,y)=u_0(s_0y),\ z(0,y)=v_0(s_0y),&0\leq y\leq 1,
\end{array}\right.
 \ees
where $f(t)$ and $g(t,y)$ are same as the proof of Theorem \ref{t1.2}. Since $u$ and $v$ are bounded, there exists a positive constant $M$ such that $|w|\leq M,\,|z|\leq M$ in $[0,+\infty)\times[0,1]$. By virtue of Lemma \ref{l3.2}, one has $0<s'(t)\leq C$ in $(0,+\ty)$ for some positive constant $C$ depending only on $s_0$, $\hat s$, $M$ and $\|u_0,\,v_0\|_{W_k^2((0,s_0))}$. So the coefficients of problem (\ref{6.7}) are bounded due to $s_0\leq s(t)<s_\infty<+\infty$. Similar to the proof of Proposition A in \cite{ZW} we can conclude that, for given $0<\sigma\ll 1$, there exists a positive constant $M_1$, which depends only on $\sigma$, $s_0$, $\hat s$, $M$ and $\|u_0,\,v_0\|_{W_k^2((0,s_0))}$, such that
  $$\|w,z\|_{C^{\frac{1+\sigma}2,1+\sigma}([1,+\infty)\times[0,1])}\leq M_1.$$
Note
 \[u_x=w_ys^{-1}(t), \ \ \ v_x=z_ys^{-1}(t), \ \ \ s'(t)=-\mu\big[(u_x(t,s(t))+\rho v_x(t,s(t))\big],\]
in view of $s(t)\geq s_0$, $0<s'(t)\leq C$ and $\|w_y(\cdot,1),z_y(\cdot,1)\|_{C^{\frac\sigma 2}([1,+\infty)}\leq M_1$, we obtain
$\|s'\|_{C^{\frac\sigma 2}([1,\infty))}\leq M_2$,
where $M_2$ depends on $\mu$, $\rho$, $C$ and $M_1$. This combined with $0<s'(t)\leq C$ suggests that $\|s\|_{C^{1+\frac\sigma 2}([1,+\infty))}\leq s_\infty+M_2$, which implies
 \[ \lim_{t\to+\infty}s'(t)=0\]
since $s'(t)>0$ and $s_\infty<\infty$.

Thanks to $u,v\geq 0$ and $\mu,\rho>0$, it is obvious that
 \[u_t-d_1u_{xx}=v^p\geq 0, \ \ v_t-d_2v_{xx}=u^q\geq 0,\ \ \ \forall \ t>0,\ \ 0<x<s(t)\]
and
 \[s'(t)>-\mu u_x(t,s(t)), \ \ s'(t)>-\mu\rho v_x(t,s(t)), \ \ \forall \ t>0.\]
One can use Lemma \ref{l6.1} to deduce (\ref{6.6}). The proof is finished.
\ \ \ \ \fbox{}

\bibliographystyle{elsarticle-num}

\end{document}